\input psfig
\def\SBIMSMark#1#2#3{
 \font\SBF=cmss10 at 10 true pt
 \font\SBI=cmssi10 at 10 true pt
 \setbox0=\hbox{\SBF Stony Brook IMS Preprint \##1}
 \setbox2=\hbox to \wd0{\hfil \SBI #2}
 \setbox4=\hbox to \wd0{\hfil \SBI #3}
 \setbox6=\hbox to \wd0{\hss
             \vbox{\hsize=\wd0 \parskip=0pt \baselineskip=10 true pt
                   \copy0 \break%
                   \copy2 \break%
                   \copy4 \break}}
 \dimen0=\ht6   \advance\dimen0 by \vsize \advance\dimen0 by 8 true pt
                \advance\dimen0 by -\pagetotal
 \dimen2=\hsize \advance\dimen2 by .25 true in
%
%
  \openin2=publishd.tex
  \ifeof2\setbox0=\hbox to 0pt{}
  \else 
     \setbox0=\hbox to 3.1 true in{
                \vbox to \ht6{\hsize=3 true in \parskip=0pt  \noindent  
                \input publishd.tex 
                \vfill}}
  \fi
  \closein2
  \ht0=0pt \dp0=0pt
 \ht6=0pt \dp6=0pt
 \setbox8=\vbox to \dimen0{\vfill \hbox to \dimen2{\copy0 \hss \copy6}}
 \ht8=0pt \dp8=0pt \wd8=0pt
 \copy8
 \message{*** Stony Brook IMS Preprint #1, #2 ***}
}

\documentstyle{amsppt}
\magnification=1200
\pagewidth{6.5 true in} \pageheight{9 true in}
\tolerance=3000
\nologo

\SBIMSMark{1996/3}{April 1996}{}
\bigskip
\bigskip
\bigskip
\bigskip
\topmatter
\title
The Periodic Points of Renormalization 
\endtitle
\author
 Marco Martens 
\endauthor
\thanks
 {Institute for Mathematical Sciences, SUNY at 
Stony Brook, Stony Brook, NY 11794-3651.} 
\endthanks
\date{April 28, 1995}\enddate
\endtopmatter
\bigskip
\bigskip
\bigskip
\bigskip
\bigskip
\bigskip
\centerline {\bf Abstract.} 

\bigskip

\flushpar 
It will be shown that the renormalization operator, acting on the space of smooth
unimodal maps with critical exponent $\alpha>1$, has periodic points of any
combinatorial type.

\vskip .5truein

\bigskip
\centerline{\bf 1. Introduction}
\bigskip

\bigskip

\flushpar
A central question in the theory of dynamical systems is whether
small scale geometrical properties of dynamical systems are universal or not,
whether they are imposed by the combinatorial properties of the systems.
The empirical discovery of universality of geometry in dynamical systems was
made by Coullet-Tresser and Feigenbaum. They studied infinitely renormalizable
period doubling unimodal maps and observed that the geometry of the invariant
Cantor set of such maps converges when looking at smaller and smaller scales.
Furthermore they observed that the limiting geometry was universal, in the 
sense that the small scale geometry of these Cantor sets  
depends only on the local behavior of the map around the critical point. 
This local behavior is specified by the critical exponent. 

\flushpar
To explain the universality of geometry they defined the period doubling 
renormalization operator on a suitable space of unimodal maps. This operator
acts like a microscope: the image under the renormalization operator is a 
unimodal map describing the geometry and dynamics on a smaller scale.
The universality of geometry could be understood by conjecturing that the
renormalization operator has a hyperbolic fixed point. In particular, the 
infinitely renormalizable unimodal maps form the stable manifold of the fixed
point of the renormalization operator. The first step in proving these 
Conjectures is showing the existence of a renormalization fixed point.

\flushpar
In 1981 Lanford showed the existence of a period doubling 
renormalization fixed point with
critical exponent $\alpha=2$. This fixed point was obtained by rigorous 
numerical analysis of the renormalization operator.
In 1986,88 H.Epstein obtained fixed points for the period doubling 
renormalization operator
with critical exponents $\alpha>1$. 

\flushpar
In 1992 D.Sullivan gave his conceptual explanation for
universality. This beautiful combination of real and complex analysis opened 
new directions in the Theory of dynamical systems. It explains the existence
of periodic points of the renormalization operator in the class of unimodal 
maps with an even critical exponent. Moreover it showed the convergence of 
renormalization.
Left was to show that these periodic points are hyperbolic. The recent work of 
Lyubich shows that there is a $1-$dimensional unstable manifold.
Furthermore McMullen showed that the infinitely renormalizable maps
form the codimension $1$ stable manifolds.
 
\flushpar
The Theory of Sullivan considers holomorphic extensions of the system on the
interval. These extensions only exist if the critical exponent is even. 
However the experiments indicate also universality for any 
other critical exponent. What has been lacking to the present is a proof of
the universality Conjectures which does not leave the interval. The
specific knowledge of the critical exponent should be irrelevant for the 
understanding of universality. The reason for wishing more general 
universality theorems is that universality is observed in many different
fields but not explained. For example universality 
in statistical mechanics, is observed but not explained. Maybe an
explanation of universality in interval dynamics which uses less structure
can spread some light on the other universality phenomena.   

\flushpar
In this work an approach to interval dynamics is developed based on real 
methods. 
The main result is the Existence Theorem of Periodic Points for the 
renormalization operator.   

\proclaim{Theorem} The renormalization operator, acting on the space of smooth
unimodal maps with critical exponent $\alpha>1$, has periodic points of any
combinatorial type.
\endproclaim

\flushpar
The Theorem is formulated precisely in section 2. This introduction will be used 
to outline the method. It will discuss the reason for the definitions given in the 
sections 3,4,...,8. The technical Lemmas are proved in section 9 and the appendix.

\bigskip
\centerline{\bf Outline of the Argument}
\bigskip

\flushpar
Let us discuss in some detail the period doubling renormalization
operator and outline the construction of a fixed point. There is no 
essential difference 
between the construction of this fixed point and the construction of a general
periodic orbit. 

\flushpar
A {\it unimodal map} is 
an endomorphism of the interval $[-1,1]$. It is of the form $f=\phi\circ q_t$
where $\phi\in \text{Diff}^2_+([-1,1])$ is an orientation preserving 
$C^2$ diffeomorphism of $[-1,1]$ and 
$q_t: [-1,1]\to [-1,1]$, $t\in [0,1]$ defined by
$$
q_t(x)=-2t|x|^\alpha+2t-1.
$$
The exponent $\alpha>1$ is called the {\it critical exponent} of $f$. 
We will fix it once and for all throughout the text. The map $q_t$ is called
the  {\it canonical folding map} with {\it peak-value} $t\in [0,1]$. The 
peak-value determines the maximum $q_t(0)=2t-1$. 
The above form for the canonical folding map is not just a 
choice for convenience. The canonical folding maps are naturally presented to 
us in section 4. They have a property intrinsically related to universality. 
The diffeomorphism $\phi$ is called the {\it diffeomorphic part} of $f$.
Observe that $f(1)=f(-1)=-1$. The collection of unimodal maps with the 
chosen critical exponent $\alpha>1$ is denoted by $\Cal{U}$.

\bigskip

\flushpar
Let $\Cal{U}^+$ be the collection of unimodal maps whose peak-value is high 
enough such that the unimodal map has a fixed point $p\in (0,1)$. For 
every $f\in \Cal{U}^+$ we can consider the first return map to the interval
$[-p,p]$. If the peak-value is not too high this first return map will be just
$f^2|_{[-p,p]}$, the unimodal map $f$ is called {\it renormalizable}. The 
unimodal 
map obtained by rescaling this first return-map to $[-p,p]$ is called the 
{\it renormalization} of $f$. The period doubling renormalization operator is 
illustrated in Figure 1. 

\midinsert
\centerline{\psfig{figure=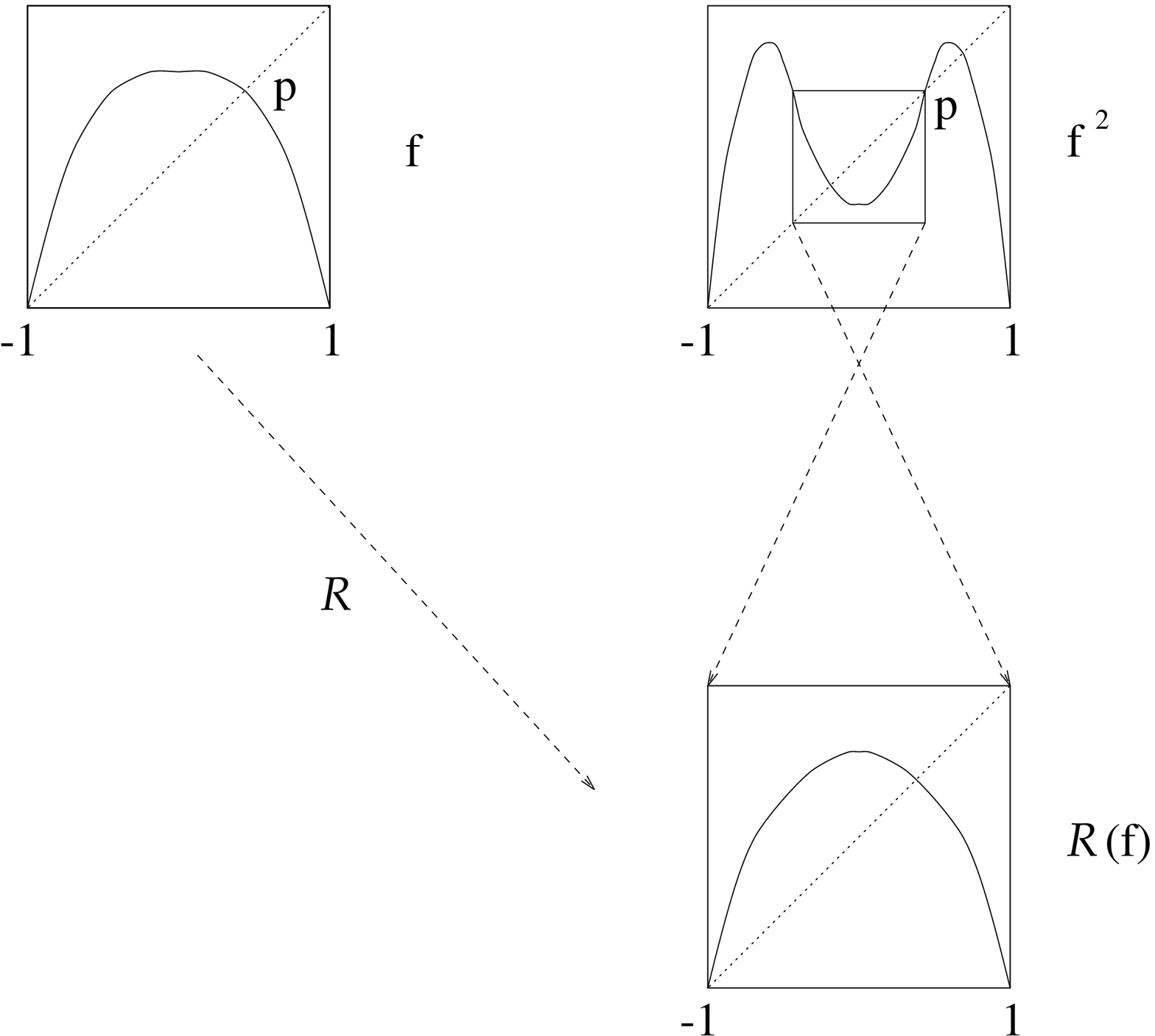,width=0.45\hsize }}
\centerline{Figure 1. The Renormalization Operator}
\endinsert

\flushpar
The Theorem states that the renormalization operator has a fixed point. The 
naive reason for the existence is the following. The space of unimodal maps
can be represented as $\Cal{U}=\text{Diff}^2_+([-1,1])\times [0,1]$. 
The renormalization operator is only defined on the subset of renormalizable 
unimodal maps. For every diffeomorphism $\phi$ there will be a range of peak-values
$t$ such that $f=\phi\circ q_t$ is renormalizable. It seems that the subset 
of renormalizable maps forms a strip, see Figure 2. This is naive but let us 
assume this situation. Maps close to the lower boundary will have 
renormalizations with peak-value close to $0$. 
Moving towards the upper boundary we will
see renormalizations whose peak-value tend to $1$. The global action of
the renormalization operator is illustrated in Figure 2.

\midinsert
\centerline{\psfig{figure=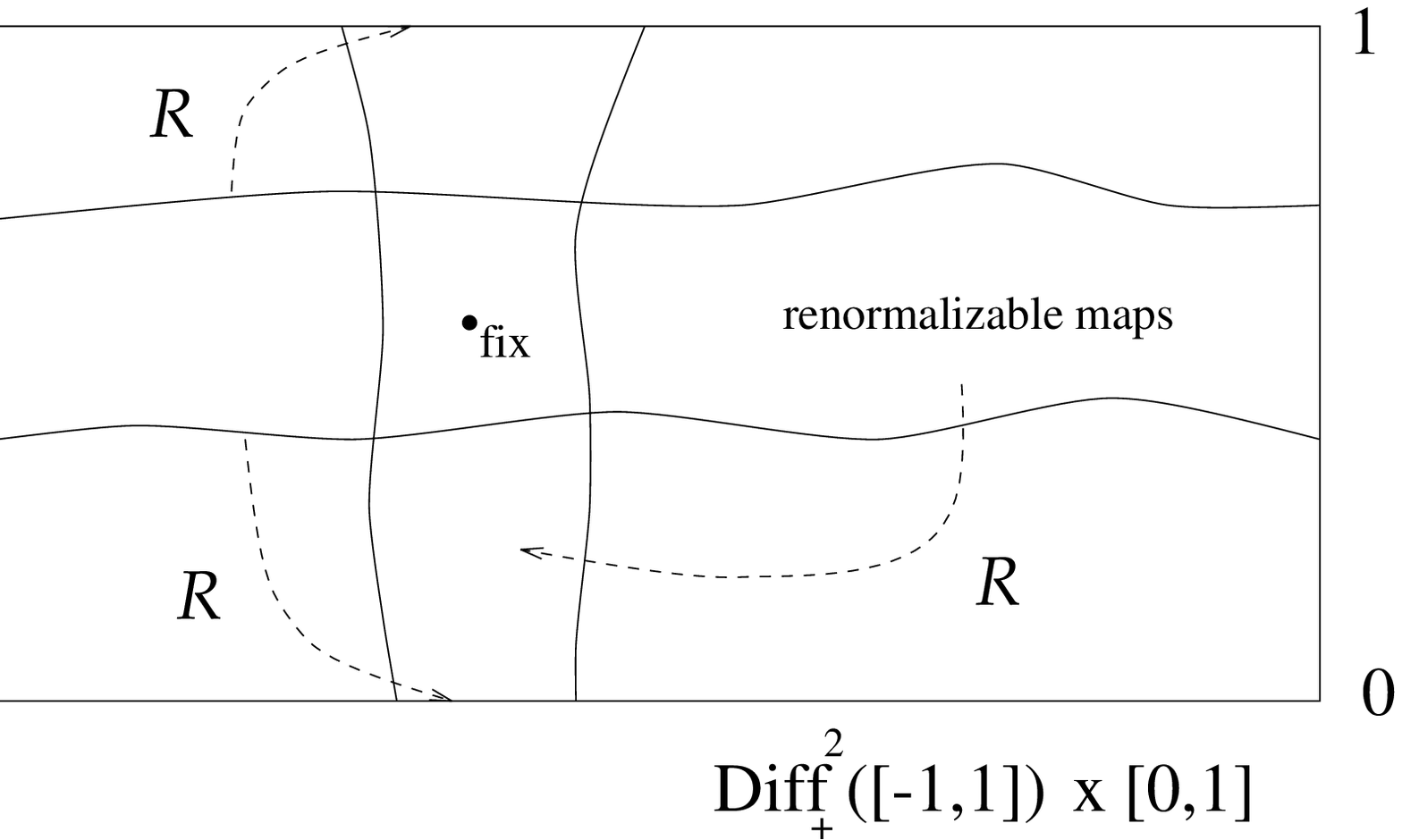,width=0.6\hsize }}
\centerline{Figure 2. Bottom Moves Down, Top Moves Up}
\endinsert

\flushpar
If we furthermore allow us to think about $\text{Diff}^2_+([-1,1])$ as being a 
compact Euclidean ball then the renormalization operator can be modeled by a map
$$
R: D\times [0,1] \to D\times (-\infty,\infty)
$$  
where $D$ is an Euclidean ball and $R(D\times \{0\})\subset D\times (-\infty,0)$
and
$R(D\times \{1\})\subset D\times (1,\infty)$: bottom goes down and top goes 
up. A slight variation on Brouwer's Fixed Point Theorem gives a fixed point 
for any map of the above form.

\flushpar
There are two difficulties. The renormalizable maps do, maybe, not form a strip
and secondly the space of diffeomorphisms is far from an Euclidean ball. There
will be a very simple way around the first problem. The way to deal with the
second problem is more elaborate and we will concentrate on this. The idea
is to construct a (thin) subspace in the space of diffeomorphism which is 
parametrized by the Hilbert-cube, a countable product of closed intervals. The
Hilbert-cube is close enough to an Euclidean ball to apply the above  
idea. This Hilbert-cube in the space of diffeomorphism is constructed in a 
natural way: it gives rise to an invariant set of the 
renormalization operator. Moreover it attracts exponentially every orbit of
the renormalization operator, the fixed point of renormalization has to be 
in this subspace. Moreover, the universality phenomena are caused by the very
special nature of these diffeomorphisms in this Hilbert-cube.
 
\bigskip

\flushpar
The construction of this attracting set of unimodal maps is based on a
generalization of the renormalization operation. The usual renormalization 
operator will, from now on,  be called the {\it classical} renormalization 
operator. A careful study of the classical renormalization operator will lead us 
to this generalization.

\flushpar 
Fix $f=\phi\circ q_t \in \Cal{U}^+$ with fixed point $p\in (0,1)$. We will 
need 
the following objects to construct the renormalization of $f$: the {\it central
interval} $S_2=[p',p]$ where 
$p'=-p$ and the {\it side interval} $S_1=[p,b]$ where $b\in (p,1)$ is such that 
$f(b)=p'$. Furthermore let $S_2(1.)=\phi^{-1}(S_2)$ and $S_1(1.)=\phi^{-1}(S_1)$ be 
the corresponding pull-backs. The notation looks too complicated. However, it
corresponds to the notation used in the general case.

\flushpar
Usually, unimodal maps are represented by their graph. The {\it dynamical
picture} in Figure 3 is more convenient for the renormalization discussion. It
contains all the objects needed to construct the renormalization.

\midinsert
\centerline{\psfig{figure=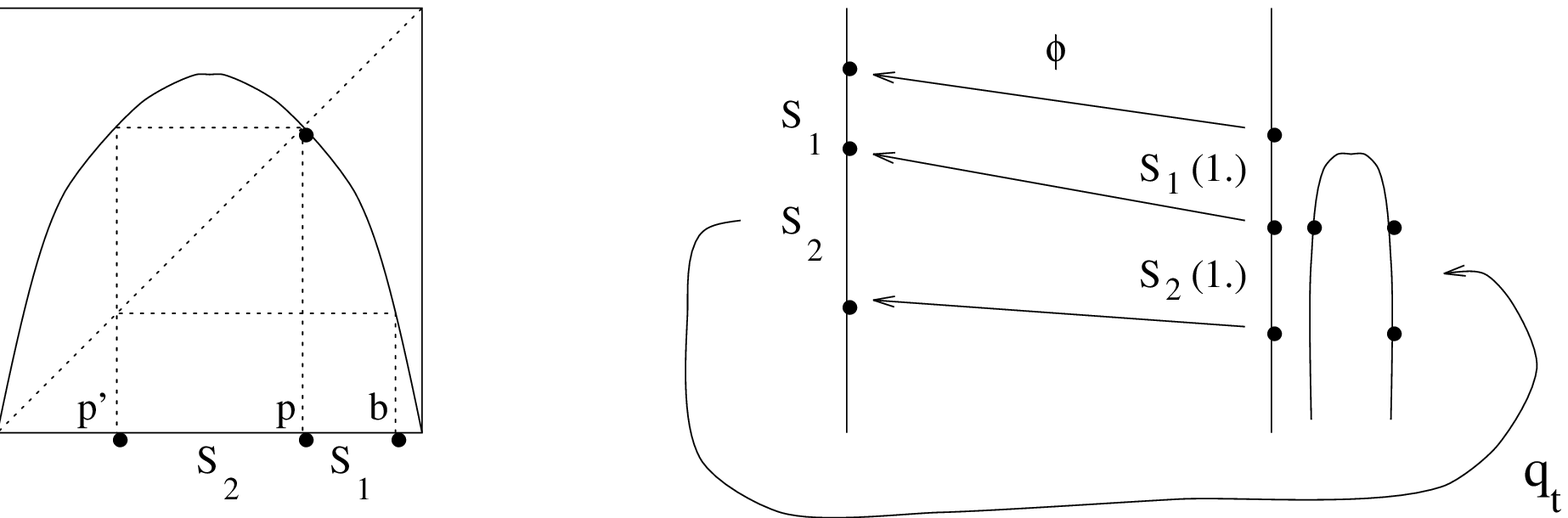,width=1.0\hsize }}
\centerline{Figure 3. A Renormalizable Unimodal Map}
\endinsert

\flushpar
The unimodal map $f$ is renormalizable if $q_t(0)\in S_1(1.)$. To describe 
precisely the classical renormalization operator we need to define the
{\it Zoom-}operators. Let $T\subset \Bbb{R}$ be a closed interval and 
$\psi:T\to \Bbb{R}$ be a diffeomorphism onto its image. Then 
$[\psi |_T]\in \text{Diff}_+^2([-1,1])$ is the orientation preserving 
diffeomorphism of $[-1,1]$ obtained by rescaling the map $\psi:T\to \psi(T)$.
If $T\subset (-1,1)$ then the operator 
$Z_T:\text{Diff}^2_+([-1,1])\to \text{Diff}^2_+([-1,1])$ defined
by 
$$
Z_T: \phi\to [\phi|_T]
$$
is called the {\it zoom-operator} to $T$.

\bigskip

\flushpar
The classical renormalization of $f$, the rescaled first return map to $S_2$, 
is given by
$$
\aligned
\Cal{R}_{\text{class}}(f)&=([\phi |_{ S_2(1.)}\circ
                              q_t |_{S_1}   \circ
                             \phi |_{S_1(1.)}] )\circ
                          q_\rho\\
                         &=([\phi |_{ S_2(1.)}]\circ
                          [q_t |_{S_1}]   \circ
                          [\phi |_{S_1(1.)}] )\circ
                          q_\rho
\endaligned
$$  
where $q_\rho: [-1,1]\to [-1,1]$ is the canonical folding part of the 
renormalization, 
obtained by rescaling domain and range of $q_t: S_2\to S_1(1.)$. Indeed this 
rescaled map is again a canonical folding map. The peak-value $\rho\in [0,1]$ 
of the renormalization is, in general, not the original peak-value $t$. 
Using the zoom-operators we get the following expression for the classical 
renormalization operator.
$$
\Cal{R}_{\text{class}}(f)=Z_{S_2(1.)}(\phi)\circ 
                          Z_{S_1}(q_t)\circ 
                          Z_{S_1(1.)}(\phi)\circ q_\rho.
$$

\bigskip

\flushpar
The classical renormalization of a unimodal map involves three operations:
\parindent=15pt
\item{1)} Finding the dynamical intervals $S_1(1.),S_2(1.)$, $S_2$ and $S_1$,  
\item{2)} Zooming in onto these dynamical intervals to obtain  
the diffeomorphisms $Z_{S_1(1.)}(\phi)$,  $Z_{S_1}(q_t)$, $Z_{S_2(1.)}(\phi)$
and the canonical folding map $q_\rho$, 
\item{3)} Compositions of these diffeomorphisms and the canonical folding map 
$q_\rho$.

\bigskip

\flushpar
The first aspect of renormalization deals with the mystery of universality.
Universality says that these intervals have very special positions. In the
construction of the fixed point the Brouwer Theorem will take care of this 
mystery, it will allow us to avoid a careful discussion. We will concentrate 
on zooming and composition.

\flushpar
The zoom operators are expected to behave well. After all, if $T\subset (-1,1)$
is a very small interval then $Z_T(\phi)$ is very close to the identity map
of $[-1,1]$. In section 2. we will define a vector space structure and the 
{\it non-linearity}-norm on $\text{Diff}^2_+([-1,1])$ which make it into the 
Banach space $\Cal{D}$. The zoom-operators on $\Cal{D}$ become linear
contractions. The second aspect of renormalization, the zooming part,
is under control.

\flushpar
The composition operator is known to behave badly. Here we propose a simple 
way around this problem: do not compose! Consider a unimodal map which can be
renormalized repeatedly. After taking a few renormalizations we will obtain a 
unimodal map whose diffeomorphic part is a long composition of interval
diffeomorphisms, restrictions of the original unimodal map. The actual 
composition of these diffeomorphisms will eliminate the information contained 
in
the factorization. We do not want to lose this information. 
This leads us to the idea of a {\it decomposition}. Instead of considering 
the diffeomorphic part of a renormalization as an interval diffeomorphism we 
will consider it to be a non-composed chain of diffeomorphisms, called a
decomposition. 

\flushpar
These decompositions are not arbitrary chains of diffeomorphisms. There is
essentially one way to label the diffeomorphisms in these chains.
Let us discuss this labeling in some detail. 

\bigskip

\flushpar
A finite ordered set $(T^{(n)},\succ)$ is called a set of 
{\it decomposition times}
if it has the following properties
\parindent=15pt
\item{1)} $T^{(n)}=\bigcup_{i=0}^n L_i$ a pairwise disjoint union of levels $L_i$,
\item{2)} $L_0=\{1.\}$,
\item{3)} there exist order preserving bijections 
$$
A_1:T_1\to \bigcup_{i=0}^{n-1} L_i \text{ and } 
A_2:T_2\to \bigcup_{i=0}^{n-1} L_i,
$$  
where $T_1=\{\tau\in T^{(n)}| \tau\prec 1.\}$ and  
      $T_2=\{\tau\in T^{(n)}| \tau\succ 1.\}$.

\flushpar
Clearly such sets exists and can be modeled by the vertices of a finite binary 
tree, see Figure 4.
A decomposition $\underline{\phi}$ is a chain of diffeomorphisms labeled 
by a set of decomposition times $T^{(n)}$, with $n\ge 0$,
$$
\underline{\phi}:T^{(n)}\to\Cal{D}.
$$
To understand the reason for this labeling we will analize how this 
labeling behaves under renormalization. We can compose the diffeomorphisms 
in the decomposition according to the order of times. Doing so we 
obtain a diffeomorphism $O(\underline{\phi})\in\Cal{D}$. Now choose a 
peak-value $t\in [0,1]$ such that the unimodal map 
$f=O(\underline{\phi})\circ q_t$ is renormalizable, see Figure 4.
We will use the notation $f=(\underline{\phi},t)$ to indicate that the
diffeomorphic part of $f$ comes from a decomposition.

\midinsert
\centerline{\psfig{figure=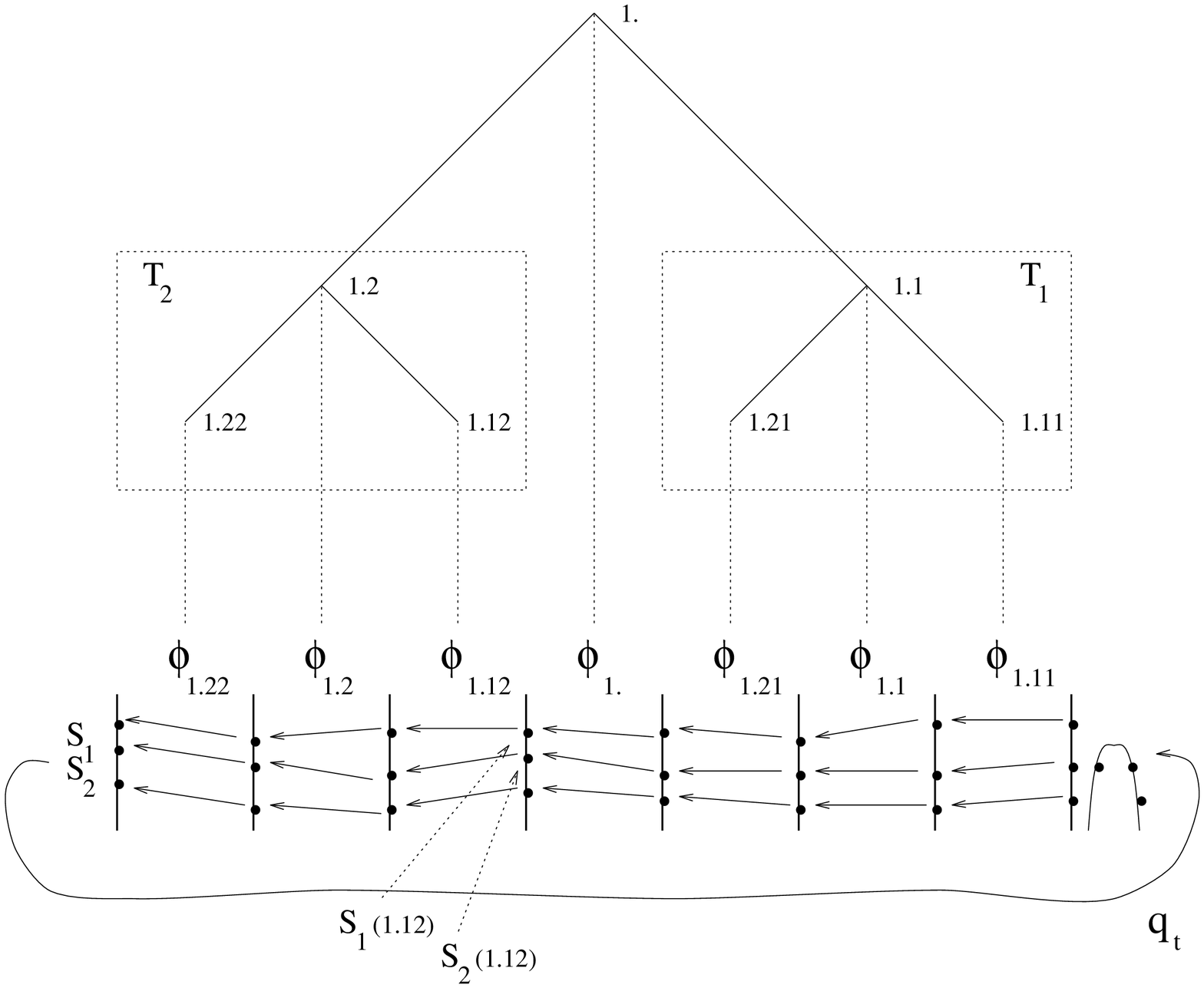,width=1.0\hsize }}
\centerline{Figure 4. The Decomposed Unimodal Map $f$}
\endinsert

\flushpar
The fixed point of $f$ is $p\in (0,1)$, the central interval $S_2=[-p,p]$ and 
the side interval $S_1=[p,b]$ are defined as before. We have to pull back 
these 
intervals through the decomposition to construct the renormalization of $f$.
These preimages of $S_2$ and $S_1$ are the dynamical intervals
$$
S_1(\tau)=(O^\tau(\underline{\phi}))^{-1}(S_1) \text{ and }
S_2(\tau)=(O^\tau(\underline{\phi}))^{-1}(S_2)
$$
where $\tau \in T^{(n)}$ and $O^\tau(\underline{\phi})\in \Cal{D}$ the 
diffeomorphism obtained by composition of the diffeomorphisms $\phi_w$
with $w\succcurlyeq \tau$, see Figure 4. 

\flushpar
The next step in the renormalization process, is to zoom in onto these 
dynamical intervals. Let 
$\underline{S}_1=(S_1(\tau))_{\tau\in T^{(n)}}$ and 
$\underline{S}_2=(S_2(\tau))_{\tau\in T^{(n)}}$. This defines the decompositions
$Z_{\underline{S}_1}(\underline{\phi}), 
Z_{\underline{S}_2}(\underline{\phi}): T^{(n)}\to \Cal{D}$ with
$$
Z_{\underline{S}_1}(\underline{\phi})(\tau)=Z_{S_1(\tau)}(\phi_\tau)
\text{ and }
Z_{\underline{S}_2}(\underline{\phi})(\tau)=Z_{S_2(\tau)}(\phi_\tau)
$$
for $\tau \in T^{(n)}$. Then form the decomposition
$$
(Z_{\underline{S}_2}(\underline{\phi}),
             Z_{S_1}(q_t), 
             Z_{\underline{S}_1}(\underline{\phi})):T^{(n+1)}\to \Cal{D},
$$
as illustrated in Figure 5.
The renormalization of $f=(\underline{\phi},t)$ becomes
$$
\Cal{R}(f)=((Z_{\underline{S}_2}(\underline{\phi}),
             Z_{S_1}(q_t), 
             Z_{\underline{S}_1}(\underline{\phi})),\rho)
$$
where $\rho$ is, as before, the peak-value of the renormalization. Compare
this with the expression for the classical renormalization. The renormalization
operation in terms of decompositions is illustrated in Figure 4. and 5.

\midinsert
\centerline{\psfig{figure=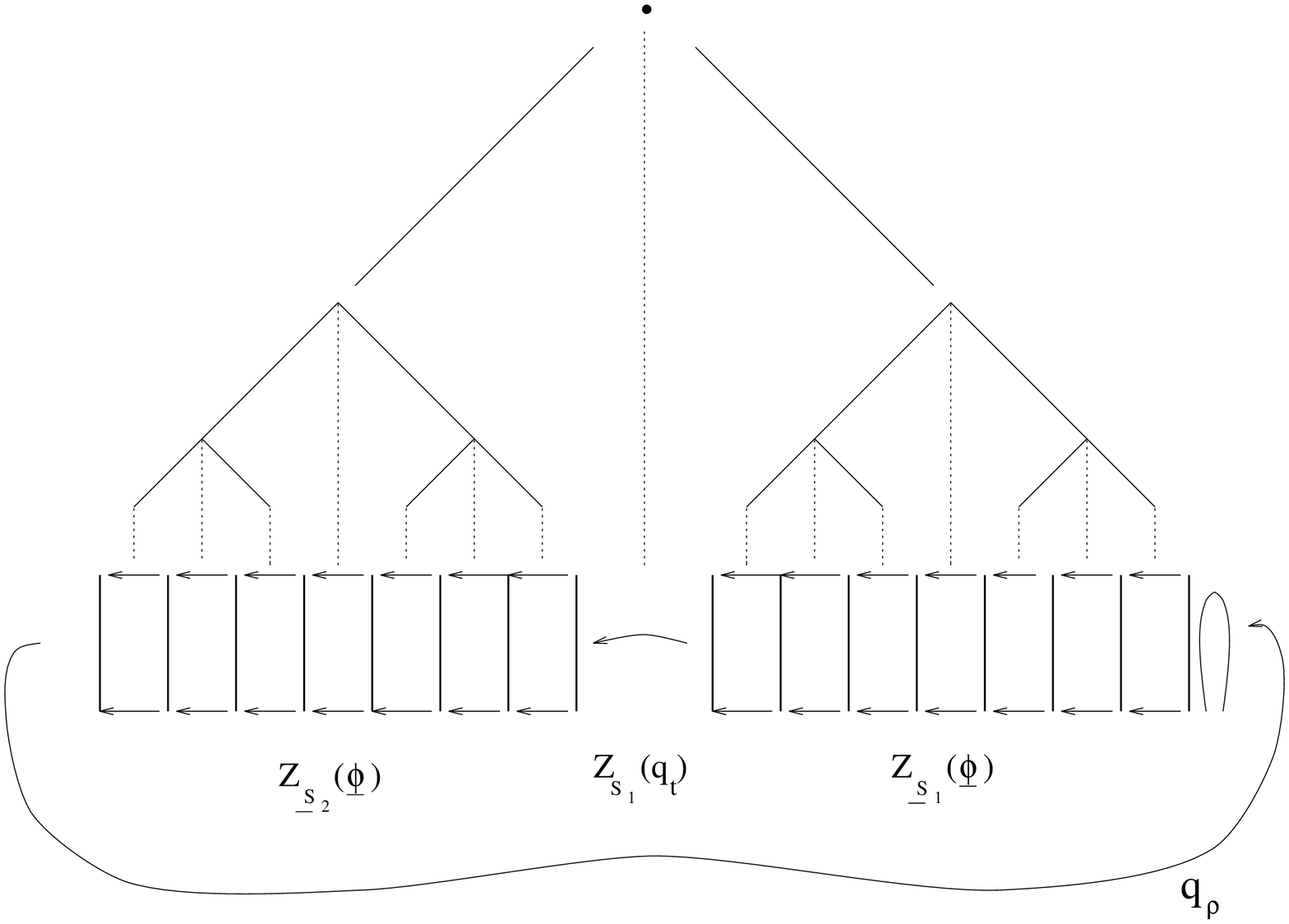,width=1.0\hsize }}
\centerline{Figure 5. The Renormalization of the Decomposed Map $f$}
\endinsert

\flushpar
Observe that the diffeomorphic part of the renormalization $\Cal{R}(f)$ is
a decomposition labeled by a tree with one level more. Each renormalization
step will give rise to decompositions labeled by a tree with one more level.
We would like to have the renormalization operator to act on a space. The
obvious definition of a decomposition will be a chain of diffeomorphisms 
labeled by the infinite binary tree $T$ or more precisely by an infinite
set of decomposition times. 

\flushpar
The Banach space of decompositions $X$, the corresponding space of 
{\it decomposed unimodal maps} $U=X\times [0,1]$ and the action of the 
corresponding renormalization operator will be discussed in the sections 4,5 
and 6.  This renormalization operator is called the 
{\it dynamical} renormalization operator $\Cal{R}_{\text{dyn}}$. 
The notion of decompositions eliminated the difficulty of the third aspect 
of renormalization, composition. 

\bigskip

\flushpar
Instead of discussing the classical renormalization operator 
$\Cal{R}_{\text{class}}$ we will discuss the dynamical renormalization 
operator $\Cal{R}_{\text{dyn}}$. This dynamical renormalization operator 
is a lift of the classical renormalization operator to the space of 
decomposed unimodal maps. In particular, the construction of a fixed point
for the dynamical renormalization operator will give a fixed point for the 
classical renormalization operator. The action of the dynamical renormalization 
operator was informally discussed above:
zoom in to dynamically defined intervals and move the obtained 
diffeomorphisms to the right place in the tree. The zoom-operators are 
contractions. If the dynamically defined intervals were always at the same
place  then the renormalization operator would essentially be a zoom-operator.
In particular renormalization would be a contraction and the 
fixed point would be found immediately.

\bigskip

\flushpar
This discussion leads to the definition of {\it geometrical renormalization
operators}. Eliminate the difficulty that the intervals
$\underline{S}_2, S_1$ and $\underline{S}_1$ have dynamical definitions. 
Make an arbitrary  choice for these 
intervals. Such a 
choice $g=(\underline{S}_2, S_1, \underline{S}_1)$ is called a {\it geometry}. The
space of all a priori possible geometries is denoted by $G$. These geometries 
are abstract objects. They have no dynamical meaning. For every geometry 
$g\in G$ define the geometrical renormalization operator
$$
\Cal{R}_g:X\to X.
$$ 
These operators act on the space of decompositions like the dynamical 
renormalization operator. The difference is that they do not zoom in onto 
dynamically defined intervals but they just zoom in onto the chosen 
intervals of the geometry $g$. They are essentially zoom-operators. 
Putting all the definitions together we see that these geometrical 
renormalization 
operators are affine contractions on the Banach space of decompositions. 
Therefore each geometrical renormalization operator has a unique fixed point. The
decompositions which are fixed point of some geometrical renormalization
operator are called {\it pure decompositions}. 
For each geometry there is a unique pure 
decomposition. These pure decompositions form a thin subspace in the space $X$ 
of decompositions, parametrized by the Hilbert-cube $G$ of geometries. 

\bigskip

\flushpar 
The geometrical renormalization operators are generalizations of the 
dynamical renormalization. Now we will explain the relation between
geometrical and dynamical renormalization and show that fixed points
of the dynamical renormalization operator have to be in the set of unimodal
maps formed by these pure decompositions. 

\flushpar
A decomposed unimodal map is a pair $f=(\underline{\phi},t)$,
where $\underline{\phi}\in X$ is a decomposition and $t\in [0,1]$ a 
peak-value. In section 4 it will be shown that decompositions can be 
composed to an actual diffeomorphism, $O(\underline{\phi})$. This 
diffeomorphism can be composed with the the canonical folding map $q_t$ and we 
obtain a classical unimodal map. It could be called the {\it observed} unimodal
map.
If the peak-value is high enough we can 
define, as before, the central interval $S_2$ and the side
interval $S_1$ depending on the decomposed unimodal map 
$f=(\underline{\phi},t)$. The collection 
$U^+\subset X\times [0,1]$ consists of those decomposed unimodal maps 
for which the central and side interval is defined.

\flushpar
For each unimodal maps in $f=(\underline{\phi},t)\in U^+$ we can pull back
the central and side interval 
through the decomposition to obtain
the dynamical intervals 
$$
\underline{S}_1(\underline{\phi},t))
\text{ and } 
\underline{S}_2(\underline{\phi},t)).
$$
This defines the {\it dynamical} geometry of $f=(\underline{\phi},t)$:
$$
d:U^+\to G
$$
with
$$
d(\underline{\phi},t)=(\underline{S}_2(\underline{\phi},t)), 
                                   S_1(\underline{\phi},t)
                     , \underline{S}_1(\underline{\phi},t))).
$$
The dynamical renormalization operator
$\Cal{R}_{\text{dyn}}:U^+\to X\times \Bbb{R}$ becomes
$$
\Cal{R}_{\text{dyn}}(\underline{\phi},t)=
(\Cal{R}_{d(\underline{\phi},t)}(\underline{\phi}),\rho(\underline{\phi},t)),
$$
where, as before,
$$
\rho:U^+\to \Bbb{R}
$$
is the peak-value of the renormalization and 
$\Cal{R}_{d(\underline{\phi},t)}$ the geometrical renormalization operator 
using the dynamical geometry $d(\underline{\phi},t)$.

\flushpar
This indicates the candidates for a renormalization fixed point. Namely if 
$(\underline{\phi},t)\in U^+$ is a renormalization fixed point then its 
diffeomorphic part $\underline{\phi}$ is 
the unique fixed point of the geometrical renormalization operator 
corresponding to the dynamical geometry $d(\underline{\phi},t)$. Decompositions 
which are the fixed point of some geometrical renormalization operator are 
called {\it pure} decompositions. Let $P$
be the space of pure decompositions. Observe that $P$ is homeomorphic to $G$.

\flushpar
Our search for a renormalization fixed point can be limited to the subspace
$$
U^+_P\subset U^+
$$
of unimodal maps whose diffeomorphic part is a pure decomposition. This space 
$U^+_P$ is homeomorphic to $P\times [0,1]$. 
A renormalization fixed point corresponds to a fixed point of
$$
\Cal{R}_{\text{dyn}}: U^+_P\approx P\times [0,1]\to P\times \Bbb{R}.
$$

\flushpar
It is the above map $\Cal{R}_{\text{dyn}}|_{U^+_P}$, 
to which we are going to apply the
Brouwer Theorem mentioned in the naive discussion in the beginning: any map
from $P\times [0,1]$ into $P\times \Bbb{R}$ whose bottom goes down and top
goes up has a fixed point. Unfortunately the space $P$ of pure
decompositions is homeomorphic to the space of geometries $G$, which  is
not exactly a Hilbert cube, it is a countable product of open intervals.

\bigskip

\flushpar
The last step in the construction is to find a priori bounds on the geometry:
If we know that 
$$
d(U_P^+)\subset G
$$ 
is compact we could apply Brouwer to this compact piece and the fixed point 
would be constructed. It turned out to be difficult to find such an
invariant piece of $P$. It is not known whether there is such a piece. We have
to proceed differently. We will work in two steps.

\flushpar
First, construct finite dimensional approximations of this map 
$\Cal{R}_{\text{dyn}}|_{U^+_P}$. 
These finite dimensional approximations, better truncations, allow a priori 
bounds. The Brouwer Theorem gives for every truncation a fixed point, 
called truncation fixed points. 

\flushpar
These truncation fixed points are almost fixed under the renormalization operator.
This extra information enables us to shown uniform bounds 
on the geometries of the truncation fixed points. These uniform 
bounds allow us to take a limit, the renormalization fixed point is found.  

\bigskip

\flushpar
{\bf Acknowledgement.} The author would like to thank M.Lyubich, D.Sullivan, 
F.Tangerman and C.Tresser for many useful renormalization discussions and 
J.Milnor for reading carefully the first version of the manuscript.

\bigskip
\centerline{\bf 2. The Classical Renormalization Operator}
\bigskip

\flushpar
In this section we will introduce combinatorial notions needed to describe
the classical renormalization operator and formulate the Existence Theorem
for Periodic Points of Renormalization. The combinatorial statements made
are folklore, their proof will be omitted. The critical exponent $\alpha>1$ 
will be fixed throughout the text.

\bigskip

\flushpar 
A unimodal map $f\in \Cal{U}$ is called {\it renormalizable } iff there exists
an expanding 
periodic point $p\in (-1,1)$ such that the first return map to the 
{\it central} interval $C=[-p,p]$ is a of the form $f^q:C\to C$ with 
$f^q(p)=p$ and $q\ge 2$. The first return map to $C$ will be, 
up to rescaling, a unimodal map. This unimodal map is a {\it renormalization} of 
$f$. Observe that a renormalization is completely 
determined by the periodic point $p$.

\flushpar
The combinatorial aspects of a renormalization are described by unimodal
permutations. A permutation on a finite ordered set is a {\it unimodal} 
permutation if the following holds. Embed the set monotonically into the 
real line. Draw the graph of the permutation. If this graph can 
be extended to the graph of a unimodal map then the permutation is called 
unimodal. 

\proclaim{Definition 2.1} A collection 
$\Cal{S}=\{S_1,S_2,\dots, S_{q-1}, S_q\}$ of oriented closed intervals
in $[-1,1]$ is called a cycle for the unimodal map $f$ if it has following 
properties 
\parindent=15pt
\item{1)} there is an expanding 
 periodic point $p\in (-1,1)$ with $S_q=[-|p|,|p|]$,
\item{2)} $f:S_i\to S_{i+1}$ , $i=1,2,\dots,q-1$, is monotone onto,
\item{3)} $f(S_q)\subset S_1$ with $f(p)\in \partial S_1$, the boundary of $S_1$,
\item{4)} the interiors of $S_1,S_2,\dots, S_q$ are pairwise disjoint,
$\Cal{S}$ inherits an order from $[-1,1]$,
\item{5)} the map 
$$
\sigma(\Cal{S}):S_i\to S_{i+1 \mod q}
$$ 
on $\Cal{S}$ is a unimodal permutation,
\item{6)} the orientation 
$$
o_{\Cal{S}}:\Cal{S}\to\{-1,1\}
$$
is such that $o_{\Cal{S}}(S_i)=1$ when $f^i(p)$ is the left boundary point of 
$S_i$ and $o_{\Cal{S}}(S_i)=-1$ otherwise.
\endproclaim

\flushpar
Observe that
\parindent=15pt
\item{1)} a unimodal map is renormalizable iff it has a cycle,
\item{2)} properties 4), 5) and 6) follow automatically once a unimodal map
has a periodic point with the first three properties. 
\item{3)} the orientation $o_{\Cal{S}}$ depends only on $\sigma$. We will
use the notation $o_{\sigma}$.

\bigskip

\flushpar
Let $\sigma$ be a unimodal permutation and 
$$
\Cal{U}_\sigma=
\{f\in \Cal{U}| f \text{ has a cycle } \Cal{S} \text{ with } 
                  \sigma(\Cal{S})=\sigma \}.  
$$ 
The unimodal maps in $\Cal{U}_\sigma$ are sometimes called 
$\sigma-$renormalizable to emphasize the type of renormalization under
consideration.
The renormalization operator
$$
\Cal{R}_{\text{class},\sigma}:\Cal{U}_\sigma\to \Cal{U}
$$
is defined to be the rescaled first return map to the smallest central interval
giving rise to a cycle $\Cal{S}$ with $\sigma(\Cal{S})=\sigma$.

\flushpar
These sets of renormalizable maps $\Cal{U}_\sigma$ are not empty, every 
family 
$t\to \phi\circ q_t\in \Cal{U}$ contains points in each 
$\Cal{U}_\sigma$.
Often a unimodal has different cycles, the sets $\Cal{U}_\sigma$ are not 
disjoint. However they are nested. For each $\sigma$ there exists a unique 
maximal factorization $\sigma=<\sigma_n,\dots,\sigma_2,\sigma_1>$ such that
$$
\Cal{R}_{\text{class},\sigma}=\Cal{R}_{\text{class},\sigma_n}\circ \dots
\circ \Cal{R}_{\text{class},\sigma_2}\circ \Cal{R}_{\text{class},\sigma_1}.
$$ 
A unimodal permutation $\sigma$ is called {\it prime} iff $\sigma=<\sigma>$.
Clearly each permutation in the maximal factorization is prime.
Using the prime unimodal permutations we obtain a partition of the
set of renormalizable unimodal maps and the {\it classical renormalization} 
operator becomes 
$$
\Cal{R}_{\text{class}}:\{\text{renormalizable maps}\}= 
                       \bigcup_{\text{prime } \sigma} \Cal{U}_\sigma
                       \to \Cal{U},
$$
with $\Cal{R}_{\text{class}}|\Cal{U}_\sigma=\Cal{R}_{\text{class},\sigma}$.

\proclaim{Theorem} For each choice $\sigma_n,\dots  \sigma_2,\sigma_1$
of prime unimodal permutations there exists $f\in \Cal{U}$ with
$$
\Cal{R}^n_{\text{class}}(f)=f
$$
and
$$
\Cal{R}^i_{\text{class}}(f)\in \Cal{U}_{\sigma_{i+1 \mod n}}
$$
for $i\ge 0$.
\endproclaim 

\flushpar
The Existence Theorem for Periodic Points is a direct consequence of 

\proclaim{Theorem 2.2} For every unimodal permutation $\sigma$ there exists
$f\in \Cal{U}_\sigma$ with
$$
\Cal{R}_{\text{class},\sigma}(f)=f.
$$
\endproclaim 

\flushpar
This Theorem 2.2 will be proved in the next sections. Fix the 
critical exponent $\alpha>1$ and the unimodal permutation $\sigma$.

\bigskip
\centerline{\bf 3. Zoom Operators}
\bigskip

\flushpar
Let $\Cal{D}=\text{Diff}_+^2([-1,1])$ be the $C^2$ orientation preserving 
diffeomorphisms of the interval $[-1,1]$. Consider the {\it non-linearity}
$N: \Cal{D}\to C^0([-1,1])$ with
$$
N(\phi)(x)=\frac{D^2\phi(x)}{D\phi(x)}=D\ln D\phi(x).
$$
This map is a bijection with inverse
$$
N^{-1}(\eta)(x)=2\frac{\int_{-1}^x e^{\int_{-1}^s\eta}ds}
                      {\int_{-1}^1 e^{\int_{-1}^s\eta}ds}-1.
$$
We will identify $\Cal{D}$ with $C^0([-1,1])$ and use the supremum norm of 
$C^0([-1,1])$. In this context we will speak about the {\it non-linearity norm}
on $\Cal{D}$. 
Observe that these linear and metric structures on $\Cal{D}$ are not
the usual structures on $\text{Diff}_+^2([-1,1])$. 
 An appendix is added in which some properties of this norm are discussed. The
 Sandwich Lemma 10.5 is the most important property. Usually we will denote 
 $N(\phi)$ by $\eta_{\phi}$.

 \bigskip

 \flushpar
 Let $I\subset [-1,1]$ be an oriented closed interval. Let $o(I)=\pm 1$,
 according to whether the orientation of $I$ and the natural orientation of
 $[-1,1]$ matches or not. Furthermore define 
 $$
 i_I:[-1,1]\to I
 $$
 to be the affine orientation preserving map onto $I$.
 Now we can define the zoom operator
 $$
 Z_I:\Cal{D}\to \Cal{D}
 $$
 by
 $$
 Z_I(\phi)=(i_{\phi(I)})^{-1}\circ \phi \circ i_I,
 $$
 where $\phi(I)$ and $I$ are oriented in the same direction, 
 $o(\phi(I))=o(I)$.

 \proclaim{Lemma 3.1} The zoom operator $Z_I$ is a linear 
 contraction. In particular
 $$
 \left| Z_I(\phi)-Z_I(\psi)\right| \le \frac{|I|}{2}\left|\phi-\psi \right|.
 $$
 \endproclaim

 \demo{Proof} If $\psi=\phi_2\circ \phi_1$ then we have the following chain 
 rule for non-linearities: 
 $$
 \eta_{\psi}=(\eta_{\phi_2}\circ\phi_1)\times D\phi_1 + \eta_{\phi_1}.
 $$
 So
 $$
 \eta_{Z_I(\phi)}=(\eta_{\phi} \circ i_I)\times Di_I=
 o(I)\frac{|I|}{2}\eta_{\phi}(i_I).
 $$
 And the statement follows.
 \hfill\hfill\qed $\,\,$ (Lemma 3.1)
 \enddemo

 \flushpar
In the next section we will use the Sandwich Lemma 10.5.  
It deals with $C^3$ diffeomorphisms. We will need 
 $$
 \Cal{D}_C=\{\phi\in\Cal{D}| \eta_{\phi}\in C^1([-1,1]) \text{ and }
		    \forall x \in [-1,1] \left|\eta'_{\phi_\tau}(x)\right|
		    \le C \left|\eta_{\phi_\tau}(x)\right|\}
 $$
 where $C>0$ is a big constant which will be defined in section 7, the proof
 of Theorem 2.2.

  \bigskip
  \centerline{\bf 4. The Space of Decompositions}
  \bigskip

\flushpar
The objects we are going to define will depend on the unimodal permutation
$\sigma$ which is fixed throughout the text. We will suppress subscripts 
$\sigma$. Let $q=|\sigma|$.
   
\bigskip

\flushpar
The set $T$ of {\it decomposition times} is a countable set with the 
properties
\parindent=15pt
\item{1)} $T$ carries an order $\prec$.
\item{2)} $\{1,2,\dots, q-1\}\subset T$ and is naturally ordered by $\prec$. 
\item{3)} the intervals 
$$
\aligned
    T_q&=\{\tau\in T| \tau\succ q-1\}\\
    T_i&=\{\tau\in T| i \succ \tau\succ i-1\}, i=2,3,\dots, q-1\\
    T_1&=\{\tau\in T|  1 \succ \tau \} 
\endaligned 
$$
admit order preserving bijections $A_i:T_i\to T$, $i=1,2,\dots q$. 

\flushpar
The set of decomposition times exists and is unique up to isomorphism. In the
period doubling case ( $q=2$ ) 
the vertices of the binary tree can be used to model $T$.
Observe that for each $q$ we get a different set of decomposition times. 

\flushpar
For each time $\tau\in T$ there is a unique number $n(\tau)\ge 0$, called the
{\it depth } of $\tau$ and a unique sequence $i_j\in \{1,2\dots, q\}$
with $j=1,2,\dots, n(\tau)$ such that
$$
A_{i_{n(\tau)}}\circ\dots 
               \circ A_{i_2}
               \circ A_{i_1}(\tau)\in \{1,2,\dots,q-1\}.
$$ 
The times in $T$ are organized in pairwise disjoint levels, 
$$
T=\bigcup_{n\ge 0} L_n,
$$
where $L_n=\{\tau\in T| n(\tau)=n\}$. The set $L_n$ is called 
the level of depth $n$.

 \bigskip

 \flushpar
 The {\it space of decompositions} is 
   $$
   X=\{(\phi_{\tau})_{\tau\in T}| \phi_\tau\in \Cal{D} \text{ and }
   \sum_{\tau\in T} |\phi_\tau|< \infty\}.
   $$
   The elements in $X$ are called {\it decompositions}.
   This set inherits the vector space structure of $\Cal{D}$. The norm will be
   $$
   |\underline{\phi}-\underline{\psi}|=
   \sum_{\tau\in T}|\phi_\tau-\psi_\tau|.
   $$
   This norm makes $(X,|.|)$ into a Banach space. 

 \bigskip

 \flushpar
 The second half of this section will define a natural notion of composition 
 of decompositions.

 \bigskip

 \flushpar
    On $X$ there are projections
   $\pi_k, \pi^\tau:X\to X$, with $k\ge 0$ and $\tau \in T$, defined by
   $$
   \pi_k(\underline{\phi})_\tau=\phi_\tau \text{ if } \tau\in \cup_{j\le k} L_j
   $$
   $$
   \pi_k(\underline{\phi})_\tau=id \text{ if } \tau\notin \cup_{j\le k} L_j
   $$
   and
   $$
   \pi^\tau(\underline{\phi})_t=\phi_t \text{ if } t\succcurlyeq \tau
   $$
   $$
   \pi^\tau(\underline{\phi})_t=id \text{ if } t\prec \tau.
   $$
   Let
   $X_k=\pi_k(X)$. On these subsets
   there are natural composition maps
   $$
   O_n, O^\tau_n:X_n\to \Cal{D},
   $$
   where both maps compose the diffeomorphisms of a decomposition according to 
 the order of the decomposition times. The map $O_n$ composes all maps of
 a given decomposition $\underline{\phi}$ and $O^\tau_n$ composes 
   all the maps $\phi_t$   with $t\succcurlyeq \tau$. Observe
 $$
 O^\tau_n=O_n\circ \pi^\tau.
 $$
 We will not be able to 
    extend
   these composition maps to the whole $X$. To extend we need some condition on 
   the derivative of the non-linearities. Choose $C>0$ sufficiently big. In 
   section 7, in the proof of Theorem 2.2, we will make the right choice for 
 $C>0$. Let 
   $$
   X_C=\{\underline{\phi}\in X| 
	  \forall \tau\in T \text{   }\eta_{\phi_\tau}\in \Cal{D}_C\}.
  $$

  \proclaim{Proposition 4.1} There exist continuous composition maps
  $$
  O, O^\tau:X_C\to \Cal{D},
  $$
  such that $O|X_C\cap X_n=O_n$ and  $O^\tau|X_C\cap X_n=O^\tau_n$. In particular
  restricted to any bounded set in $X_C$, all composition maps $O,O^\tau$ are
  Lipschitz with the same constant.
  \endproclaim

\demo{Proof} The proof relies havely on the Sandwich Lemma 10.5, given 
in the appendix.
Let us first show that the pointwise limit of $O_n\circ \pi_n$ is defined. 
Let $\underline{\phi}\in X_C$. We are going 
to apply the Sandwich Lemma with $b=|\underline{\phi}|$ and $C$ the defining
constant for $X_C$. 
Observe that for $k>n$,
$O_k(\pi_k(\underline{\phi}))$ is obtained from 
$O_n(\pi_n(\underline{\phi}))$ by 
applying the Sandwich Lemma to all $\phi_\tau$ with 
$\tau \in \bigcup_{j=n+1}^k L_j$: 
the sequence $O_n(\pi_n(\underline{\phi}))$ 
is a Cauchy sequence. It converges to $O(\underline{\phi})$. This defines 
$$
O:X_C\to \Cal{D}.
$$
Left is to prove that the function $O$ 
is Lipschitz on bounded sets. This will be again a
Sandwich argument. Fix a bounded set $B\subset X_C$. Take 
$\underline{\phi},\underline{\psi}\in B$. The chain-rule for 
non-linearities gives for every $\tau\in T$
$$
\aligned
\left|\psi_{\tau}\circ (\phi_\tau)^{-1}\right|&\le
\sup_{x}\left| \eta_{\psi_\tau}(\phi_{\tau}^{-1}(x))-
                            \eta_{\phi_\tau}(\phi_{\tau}^{-1}(x))\right|
                           \cdot (\phi_{\tau}^{-1})'(x)\\
&\le K\cdot |\psi_\tau-\phi_\tau|,
\endaligned
$$
where the second estimate was obtained by applying Lemma 10.3. We are 
going to apply the Sandwich Lemma 10.5 with $b=K \cdot \text{diam}(B)$ and 
$C$ the defining constant for $X_C$. 
Consider again $\underline{\phi},\underline{\psi}\in B $. 
Because $O(\pi_k(\underline{\phi}))=O_k(\pi_k(\underline{\phi}))$ 
we can find for each $\epsilon>0$
some $k\ge 0$ such that
$$
\left|O(\underline{\phi})-O(\pi_k(\underline{\phi}))\right|\le \epsilon
$$
and
$$
\left|O(\underline{\psi})-O(\pi_k(\underline{\psi}))\right|\le \epsilon.
$$
Then
$$
\left|O(\underline{\psi})-O(\underline{\phi}) \right|\le
\left|O(\underline{\psi})-O(\pi_k(\underline{\psi})) \right|+
\left|O(\pi_k(\underline{\psi}))-O(\pi_k(\underline{\phi})) \right|+
\left|O(\pi_k(\underline{\phi}))-O(\underline{\phi}) \right|.
$$
But $O(\pi_k(\underline{\psi}))$ is obtained from
$O(\pi_k(\underline{\phi}))$ by Sandwiching the maps 
$\psi_\tau\circ (\phi_\tau)^{-1}$. We get
$$
\aligned
\left|O(\underline{\psi})-O(\underline{\phi}) \right|&\le
\epsilon + \text{const}\cdot 
\sum_{\tau} \left| \psi_\tau\circ (\phi_\tau)^{-1}\right|
+\epsilon\\
&\le 2\epsilon+\text{const}\cdot
\sum_{\tau} \left|\psi_\tau-\phi_\tau\right|\\
&\le 2\epsilon + \text{const}\cdot
\left|\underline{\psi}-\underline{\phi}\right|.
\endaligned
$$
Because $\epsilon$ was taken arbitrarily we get the Lipschitz 
estimate for $O$. Let 
$$
O^\tau=O\circ \pi^\tau.
$$
Because the projection do not increase distance we get the same Lipschitz
constant for all $O^\tau$, $\tau\in T$.
\hfill\hfill\qed $\,\,$ (Proposition 4.1)
\enddemo

\bigskip
\centerline{\bf 5. The Space of Geometries}
\bigskip

\flushpar
The orientation of a cycle with $\sigma(\Cal{S})=\sigma$ 
is denoted by $o_\sigma$.

\proclaim{Definition 5.1} Let $\epsilon \ge 0$. 
A collection $\Cal{S}=\{S_1,S_2,\dots,S_q\}$ of oriented closed intervals in 
$[-1,1]$ is called an $\epsilon-$elementary geometry if it has the following 
properties.  
\parindent=15pt
\item{1)} the interiors of the intervals are pairwise disjoint, 
$\Cal{S}$ inherits an order from $[-1,1]$,
\item{2)} the permutation $S_i\to S_{i+1 \mod q}$ on the ordered set $\Cal{S}$
is isomorphic to the unimodal permutation $\sigma$,
\item{3)} the orientation $o(S_i)=o_\sigma(i)$, $i=1,2,\dots,q$,
\item{4)} $\bigcup S_i \subset [-1+\epsilon,1-\epsilon]$. 
\endproclaim

\flushpar
The space of $\epsilon-$elementary 
geometries is denoted by $E_\epsilon$. Moreover 
$$
Q_\epsilon=\{\Cal{S}\in E_\epsilon | S_q=[-|p|,|p|]
\text{ for some } p\in (-1,1) \text{ with } |p|\ge \epsilon\}.
$$
The spaces $E_\epsilon$ and $Q_\epsilon$ can be considered to be convex 
subsets of some Euclidean space. The Euclidean topology makes them into 
compact Euclidean balls. The space of {\it $\epsilon-$geometries} is 
$$
G_\epsilon=Q_\epsilon\times E_\epsilon^T. 
$$
We will use the product topology on $G_\epsilon$. In particular it is 
compact. We will use the following notation: $G=G_0$ and if 
$g=(\Cal{S},(\Cal{S}(\tau))_{\tau\in T})\in G$ then
$$
|g|=\frac12 \sup \{ 
\{ \sum_{i=1}^q \left| S_i(\tau) \right| |\tau\in T\}
\cup 
\{ \sum_{i=1}^q \left| S_i \right|\} \}.
$$

\bigskip

\flushpar
The geometries do not have any dynamical meaning. They are merely 
abstract generalizations of cycles. To explain this we go back to unimodal
maps.

\bigskip

\flushpar
The non-linearity of $x\to |x|^\alpha$ is called the 
{\it critical non-linearity}
$$
\gamma(x)=(\alpha-1)\frac{1}{x}.
$$
The {\it canonical folding family} is the family
$q_t:[-1,1]\to [-1,1]$, $t\in [0,1]$,  with the properties
\parindent=15pt
\item{1)} $\eta_{q_t}=\gamma$.
\item{2)} $q_t(-1)=q_t(1)=-1$.
\item{3)} $t=\frac{|[-1,q_t(0)]|}{|[-1,1]|}$.

\flushpar
A computation gives
$$
q_t(x)=-2t\cdot |x|^\alpha+2t-1.
$$
There are two observations to be made.
First, every $q_t$ has negative Schwarzian derivative. Furthermore
consider the interval $I=[-p,p]\subset [-1,1]$ and consider $Z_I(q_t)$.
It has non-linearity
$$
\gamma (i_I(x))\cdot p=(\alpha-1)\cdot \frac{1}{px}\cdot p=\gamma (x)
$$
The canonical folding family have the fundamental property that it is fixed
under Zoom-operators to intervals centered around $0$. 

\flushpar
Usually a unimodal map is defined as a map  $f=\phi\circ q_t$
with $t\in [0,1]$. Here $\phi$ is some orientation preserving 
$C^2$ diffeomorphism of $[-1,1]$. For renormalization purposes
it is more convenient to use the following 

\proclaim{Definition 5.2} A {\it decomposed unimodal map} is a pair 
$(\underline{\phi},t)$, where $\underline{\phi}$ is a decomposition in $X_C$
and $t\in [0,1]$. The decomposition $\underline{\phi}$ is called the 
diffeomorphic part of the decomposed unimodal map $(\underline{\phi},t)$.
The interpretation of $(\underline{\phi},t)$ is classical unimodal map
$$
f=O(\underline{\phi})\circ q_t\in \Cal{U}.
$$
The space of decomposed unimodal maps is denoted by $U=X_C\times [0,1]$, with
the product topology.
\endproclaim

\proclaim{Definition 5.3} A decomposed unimodal map $f=(\underline{\phi},t)$ 
is called quasi renormalizable if the classical unimodal map
$O(\underline{\phi})\circ q_t$ has the following properties 
\parindent=15pt
\item{1)} there exists an expanding periodic point $p\in (-1,1)$, $f^q(p)=p$,
\item{2)} there exist $\epsilon>0$ and an $\epsilon-$elementary geometry
$\Cal{S}=\{S_1,S_2,\dots, S_{q-1},S_q=[-|p|,|p|]\}\in Q_\epsilon$ such that
\item{3)} $O(\underline{\phi})\circ q_t(S_i)=S_{i+1}$, $i=1,2,\dots, q-1$.
\item{4)} $|p|>0$ is minimal with the above properties.

\flushpar
A decomposed unimodal map is called renormalizable if it is quasi
renormalizable and 
\item{5)} $O(\underline{\phi})\circ q_t(S_q)\subset S_1$.
\endproclaim

\flushpar
Observe that the elementary geometry $\Cal{S}$ is uniquely defined.
It is constructed by pulling back the central interval $S_q=[-|p|,|p|]$ according
to the unimodal permutation $\sigma$.

\bigskip

\flushpar
Let $U^+\subset U$ be the set of quasi-renormalizable decomposed 
unimodal maps. The {\it dynamical geometry}
$$
d:U^+\to G
$$
is defined as follows. Let $\Cal{S}$ be the elementary geometry of  
$f=(\underline{\phi},t)$. Then
$$
d(\underline{\phi},t)=(\Cal{S}, (\Cal{S}_\tau)_{\tau\in T}),
$$
where
$$
O^\tau(\Cal{S}_\tau)=\Cal{S}
$$
for all $ \tau\in T$. The {\it peak-value} of the renormalization
$$
\rho:U^+\to[0,\infty)
$$
is defined as follows. Let $\hat{S}_1=(O(\underline{\phi}))^{-1}(S_1)$. Then
$$
\rho(\underline{\phi},t)=\frac{|q_t(S_q)|}{|\hat{S}_1|}.
$$
A decomposed unimodal map $f=(\underline{\phi},t)\in U^+$ is
renormalizable if 
$$
q_t(0)\in \hat{S}_1,
$$ 
or equivalently
$$
\rho(\underline{\phi},t)\in [0,1].
$$

\proclaim{Lemma 5.4} The dynamical geometry and the function $\rho$
are  continuous.
\endproclaim

\demo{Proof} Let $(\underline{\phi}_n, t_n)\to (\underline{\phi},t)$. 
This implies
that $\Cal{S}(n)$ will tend to $\Cal{S}$ of the limit because of the continuity 
of the composition operator $O$. Now we have to pull back those intervals by
the partial compositions $O^\tau (\underline{\phi}_n)$, $\tau\in T$. 
Observe that we see a uniform convergence of  
$O^\tau (\underline{\phi}_n)\to O^\tau (\underline{\phi})$. Lemma 10.2 and 
Lemma 10.4 can be used to show that the uniform convergence implies that the 
elementary geometries $\Cal{S}_\tau(n)$ tend uniformly to the elementary 
geometries of the limit. This means that the 
functions $d$ and $\rho$ are continuous. 

\hfill\hfill\qed $\,\,$ (Lemma 5.4)
\enddemo

\bigskip
\centerline{\bf 6. Geometrical Renormalization}
\bigskip

\flushpar
Let $g\in G_\epsilon$, say
$$
g=(\Cal{S}, (\Cal{S}_\tau)_{\tau\in T})
$$
with $\Cal{S}=\{S_1,S_2,\dots,S_q\}$ and 
$\Cal{S}_\tau=\{S_1(\tau), S_2(\tau),\dots, S_q(\tau)\}$ and 
$$
a_i=Z_{S_i}(q_t)\in \Cal{D}, i=1,2,\dots, q-1.
$$
Observe that these diffeomorphisms $a_i$ depend only on the intervals $S_i$,
not on $t$. The {\it geometrical renormalization} operator for geometry 
$g\in G_\epsilon$ with $\epsilon>0$ 
$$
\Cal{R}_g:X\to X
$$
is defined by
$$
\aligned
\Cal{R}_g(\underline{\phi})(\tau)&=
Z_{S_1(A_1(\tau))}(\phi_{A_1(\tau)}) \text{ for } \tau\in T_1\\
   \Cal{R}_g(\underline{\phi})(1)&=a_1\\
\Cal{R}_g(\underline{\phi})(\tau)&=
Z_{S_2(A_2(\tau))}(\phi_{A_2(\tau)}) \text{ for } \tau\in T_2\\
   \Cal{R}_g(\underline{\phi})(2)&=a_2\\
                            \dots&= \dots\\
\Cal{R}_g(\underline{\phi})(q-1) &=a_{q-1}\\
\Cal{R}_g(\underline{\phi})(\tau)&=
Z_{S_q(A_q(\tau))}(\phi_{A_q(\tau)}) \text{ for } \tau\in T_q,
\endaligned
$$
This definition corresponds to the informal discussion in the introduction
applied to a general dynamical picture like the one shown in Figure 6. 

\midinsert
\centerline{\psfig{figure=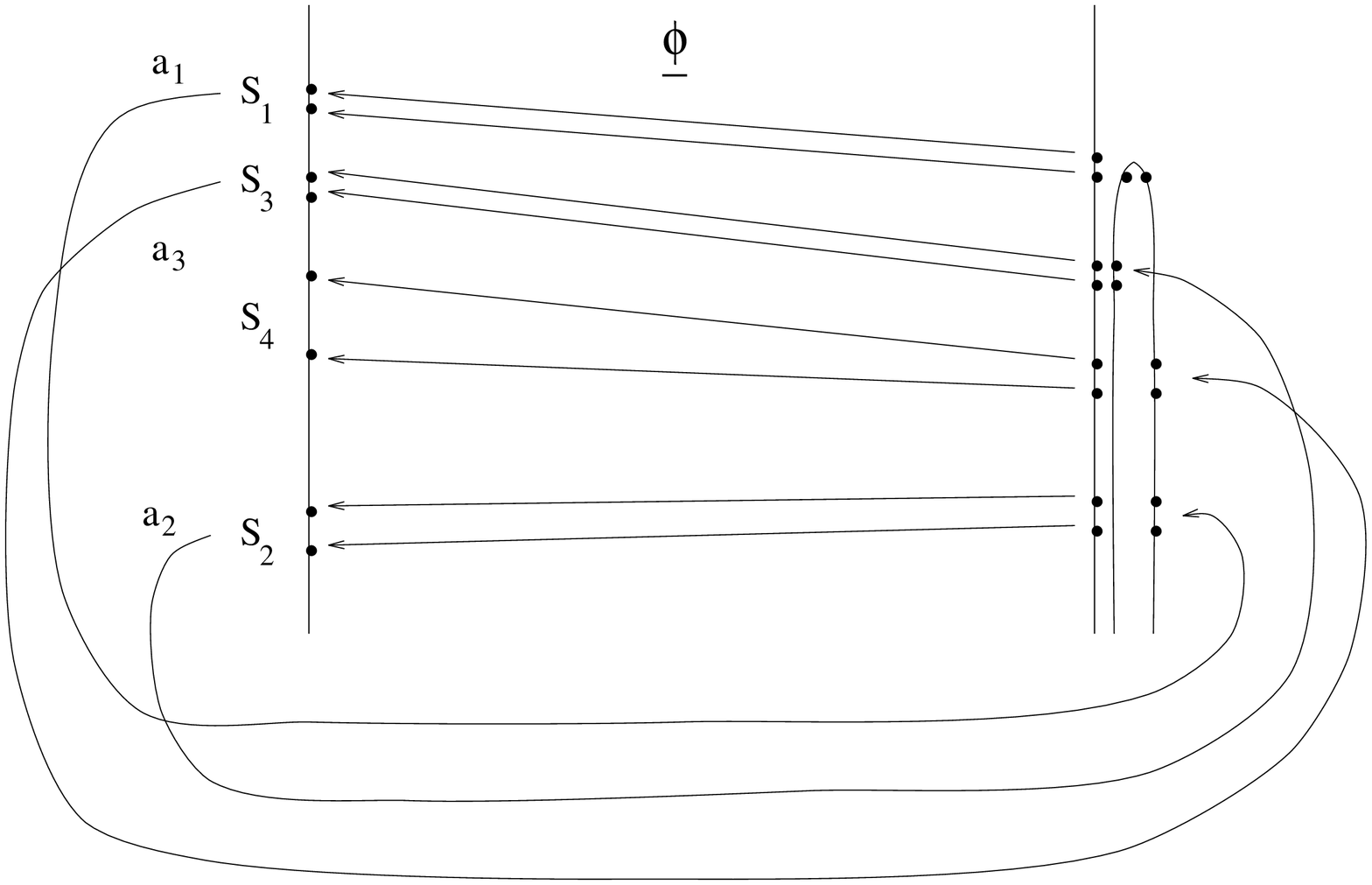,width=.75\hsize }}
\centerline{Figure 6. A General Dynamical Picture with $q=4$}
\endinsert

\bigskip

\flushpar
The family of geometrical renormalization operators is denoted by
$$
\Cal{R}:X\times G_\epsilon \to X,
$$
with $\Cal{R}(\underline{\phi},g)=\Cal{R}_g(\underline{\phi})$. 
The {\it dynamical renormalization} operator
$$
\Cal{R}_{\text{dyn}}: U^+\to X\times [0,\infty),
$$
reflects the classical renormalization operation in the space of decomposed 
unimodal maps. It is defined as follows. Let $f=(\underline{\phi},t)$ be a 
decomposed unimodal map in $U^+$ then 
$$
\Cal{R}_{\text{dyn}}((\underline{\phi},t))=
(\Cal{R}_g(\underline{\phi}),\rho(\underline{\phi},t)),
$$
where $\Cal{R}_g$ is the geometrical renormalization operator with
geometry $g=d(\underline{\phi},t)$. 

\flushpar
The three 
renormalization operators $\Cal{R}_{\text{class}}$, $\Cal{R}_{\text{dyn}}$ 
and 
$\Cal{R}$ are related by
$$
\aligned
\Cal{R}_{\text{class},\sigma }(O(\underline{\phi})\circ q_t)&=
O(\Cal{R}(\underline{\phi}, d(\underline{\phi},t)))
\circ q_{\rho(\underline{\phi},t)},\\
\Cal{R}_{\text{dyn}}(\underline{\phi},t)
          &=(\Cal{R}(\underline{\phi}, d(\underline{\phi},t)),
                      \rho(\underline{\phi},t)),
\endaligned
$$ 
where $(\underline{\phi},t)\in U^+$ is a quasi renormalizable 
decomposed unimodal map.
The difference between the dynamical and geometrical renormalization
operators is that the first uses dynamically defined geometries in the 
renormalization process and the others use given geometries.

\bigskip

\flushpar
The following Proposition is a central ingredient for the understanding of
the classical renormalization operator.

\proclaim{Proposition 6.1} For every $g\in G_\epsilon$ the geometrical
renormalization operator $\Cal{R}_g$ is an affine contraction with contraction
constant $|g|\le 1-\epsilon$. Its fixed point $\Phi(g)$ depends 
continuously on $g$. In particular $\Phi(G_\epsilon)$ is compact.
The fixed points have the following special properties
\parindent=15pt
\item{1)} Every map $O(\Phi(g))$ expands hyperbolic distance.
\item{2)} There exists $C>0$, depending on $\epsilon$, such that 
$\Phi(G_\epsilon)\in X_C$.
\endproclaim

\flushpar
The use of hyperbolic distance in interval dynamics is thoroughly discussed in
[MS]. 

\demo{Proof} The fact that $\Cal{R}_g$ is an affine contraction is a 
direct consequence of Lemma 3.1. Observe that the affine term of $\Cal{R}_g$
is formed by the diffeomorphisms $a_i$, $i=1,2,\dots, q-1$.

\flushpar
The second special
property can be shown as follows. Observe that the diffeomorphisms $a_i$, 
$i=1,2,\dots,q-1$ are obtained by applying zoom-operators to canonical 
folding maps. Their non-linearities have explicit formulas and there is some 
$C>0$ such that for all $g=(\Cal{S},(\Cal{S}_\tau)_{\tau\in T})\in G_\epsilon$
the diffeomorphisms $a_i$, 
$i=1,2,\dots,q-1$  will satisfy
$$
|\eta'_{a_i}(x)|\le C |\eta_{a_i}(x)|
$$
for all $x\in [-1,1]$. Now observe that for any oriented interval 
$I\subset [-1,1]$ the diffeomorphisms $Z_I(\eta_{a_i})$ will
satisfy the same property, because
$$
|Z_I(\eta_{a_i})(x)|=|\eta_{a_i}(i_I(x))|\cdot \frac{|I|}{2}
$$
and 
$$
|(Z_I(\eta_{a_i}))'(x)|=|\eta'_{a_i}(i_I(x))|\cdot \frac{|I|^2}{4}.
$$
It is left is to show that $\Phi$ is continuous. Observe that for any 
$k\ge 0$ 
$$
\pi_k(\Phi(g))=\Cal{R}^{k+1}_g(id)
$$
and that this projection depends only on the elementary geometries in
the levels $L_0, \dots, L_{k}$  and on $\Cal{S}$. The map 
$$
g\mapsto \pi_k(\Phi(g))
$$
is continuous. The last things to observe is that for 
$\underline{\phi}=\Phi(g)$, $g\in G_\epsilon$
$$
\sum_{\tau\in L_k} |\phi_\tau|\le (1-\epsilon)^k\cdot \sum_{i=1}^{q-1}|a_i|,
$$
which is true be cause $\Cal{R}_g$ is an $(1-\epsilon)$-contraction. So
for $\underline{\phi}=\Phi(g)$ we get
$$
\sum_{\tau\in\cup_{j\ge k} L_j} |\phi_\tau|\le \frac{(1-\epsilon)^k}
                                                    {\epsilon}\cdot 
 \sum_{i=1}^{q-1}|a_i|
\le \frac{(1-\epsilon)^k}{\epsilon^2}\cdot (\alpha-1)\cdot (q-1).
$$
To get two fixed points close, we only have to get them close on 
finite levels, which can be done because of the continuity the map
$g\mapsto \pi_k(\Phi(g))$.

\flushpar
The first special property is a consequence
of the fact that the diffeomorphisms $a_i$, $i=1,2,\dots,q-1$,
 have negative Schwarzian derivative. Observe 
that $\Phi(g)=\lim_{k\to\infty} \Cal{R}^k(id)$, all the diffeomorphisms
of the decomposition $\Phi(g)$ are obtained by applying zoom-operators
to the diffeomorphisms $a_i$, they all have 
negative Schwarzian derivative. Then observe that all finite compositions
expands definitely hyperbolic distance, which is preserved by taking the
limit.
\hfill\hfill\qed~~(Proposition~6.1)
\enddemo

\flushpar
The diffeomorphisms $O(\Phi(g))$, with $g\in G_\epsilon$, are $C^2$ by
construction. It can be shown that these diffeomorphisms are in fact
analytic.

\bigskip
\centerline{\bf 7. The fixed Point}
\bigskip

\flushpar
In this section we are going to prove Theorem 2.2.  The aim is to construct a 
fixed point for $\Cal{R}_{\text{class},\sigma}$. The strategy is to construct a
sequence of approximate fixed points for the corresponding renormalization
operator $\Cal{R}_{\text{dyn}}$ on the space of decomposed unimodal maps.
These approximate fixed points will be called truncation fixed points. 
The limit of these approximations 
is going to be a fixed point for $\Cal{R}_{\text{dyn}}$. After composition
we obtain a fixed point for $\Cal{R}_{\text{class}, \sigma}$. 

\bigskip

\flushpar
The {\it critical value} of a decomposed unimodal map is given by the continuous 
function
$$
v:\bigcup_{C<\infty} X_C\times \Bbb{R}\to \Bbb{R},
$$ 
where 
$$
\aligned
v(\underline{\phi},t)&=-1+t \text{ for } t\le 0\\
v(\underline{\phi},t)&=O(\underline{\phi})\circ q_t(0) \text{ for } t\in [0,1]\\ 
v(\underline{\phi},t)&=t \text{ for } t\ge 1.
\endaligned
$$
Note that $v$ is strictly monotone in $t$.

\flushpar
The proofs of the 
following Propositions are somewhat involved and we postpone them to 
respectively section 8. and 9. The first states the existence of truncation
fixed points and the second the a priori bounds on their dynamical geometry. 

\proclaim{Proposition 7.1} For every $k\ge 0$ there exists a decomposed 
unimodal map $(\underline{\phi},t)$ with the following properties
\parindent=15pt
\item{1)} $(\underline{\phi},t)\in U^+$,
\item{2)} $v(\underline{\phi},t)=v(\Cal{R}_{\text{dyn}}(\underline{\phi},t))$,
\item{3)} $d(\underline{\phi},t)=g\in G_\epsilon$ with $\epsilon>0$,
\item{4)} $\underline{\phi}=\Cal{R}_g^{k+1}(id)$,
\item{5)} In particular
$$
\Cal{R}_{\text{dyn}}(\underline{\phi},t)
=(\Cal{R}_g(\underline{\phi}),\rho(\underline{\phi},t))
= (\Cal{R}_g^{k+2}(id),\rho(\underline{\phi},t)).
$$

\flushpar
A decomposed unimodal map with these properties is called a truncation fixed point 
of depth $k\ge 0$.
\endproclaim 

\proclaim{Proposition 7.2} The dynamical geometry of any truncation fixed point is 
contained in a universal $G_\epsilon$, $\epsilon>0$. 
\endproclaim

\flushpar
We will use the following notation to describe truncation fixed points. For 
$k\ge 0$ define the function
$$
\Phi_k:\bigcup_{\epsilon>0} G_\epsilon \to \bigcup_{C<\infty} X_C
$$
by
$$
\Phi_k(g)=\Cal{R}_g^{k+1}(id)=\pi_k\circ \Phi(g). 
$$
These functions describe truncations of the fixed points $\Phi(g)$. Observe 
that
$\Phi_k(g)$ and $\Phi_{k+1}(g)$ differ only in the level of depth $k+1$. 

\proclaim{Theorem 2.2} The renormalization 
operator $\Cal{R}_{\text{class},\sigma}$ (and $\Cal{R}_{\text{dyn}}$) has 
a fixed point.
\endproclaim

\demo{Proof} From Proposition 7.1 we get a sequence of truncation fixed points with
increasing depth $k\ge 0$, $(\Phi_k(g_k),t_k)\in U^+$. According to Proposition 7.2
the sequence of geometries $g_k\in G_\epsilon$ can be assumed to be convergent
$$
(g_k,t_k)\to (g,t)\in G_\epsilon\times [0,1]
$$
for $k\to\infty$. This convergence implies 
\parindent=15pt 
\item{1)} $\Phi_k(g_k), \Phi_{k+1}(g_k)\to\Phi(g)$,
\item{2)} $\left| \Phi_k(g_k)-\Phi_{k+1}(g_k)\right|\le K(1-\epsilon)^k$,
where $\epsilon>0$ and $K$ are universal,
\item{3)} $\rho_k=\rho(\Phi_k(g_k),t_k)\to t$.

\flushpar
The first two statements follow from the fact that for all 
$\underline{\phi}=\Phi(g)$ with $g\in G_\epsilon$
$$
\sum_{\tau\in L_j} |\phi_\tau| \le K(1-\epsilon)^j,
$$ 
with $j\ge 0$. The third statement follows from continuity of $v$: observe
$$
v(\Phi_{k+1}(g_k),\rho_k)=v(\Phi_k(g_k),t_k)\to v(\Phi(g),t).
$$
So
$$
(\Phi_{k+1}(g_k),\rho_k)\to v^{-1}(v(\Phi(g),t))\cap(\Phi(g)\times [0,1])=
\{(\Phi(g),t)\}.
$$

\bigskip

\flushpar
The candidate fixed point for $\Cal{R}_{\text{dyn}}$ is $(\Phi(g),t)$. 
First we have to show that
$(\Phi(g),t)\in U^+$. Let $g=(\Cal{S}, (\Cal{S})_{\tau\in T})$ and 
$S_q=[-|p|,|p|]$. By continuity we see that $p$ is a periodic point of
$O(\Phi(g))\circ q_t$. If this periodic point is expanding then the elementary
geometry $\Cal{S}$ makes $(\Phi(g),t)$ renormalizable.
We have to show that $p$ is an expanding periodic point. Let
$$
B_0=v^{-1}([-1,0])\cap (\Phi(G_\epsilon)\times [0,1]).
$$
This set is compact. Furthermore let
$$
\overline{U^+_\epsilon}=\overline{U^+\cap (\Phi(G_\epsilon)\times [0,1])}
$$
This set is also compact and $B_0\cap U^+_\epsilon=\emptyset$. The distance
of these sets is $\delta>0$.

\flushpar
Assume $(\Phi(g),t)\notin U^+$. Then the periodic point $p$ is a neutral
periodic point for the map $f=O(\Phi(g))\circ q_t$. This map expands 
hyperbolic distance. So the periodic orbit attracts $0$. The orbit of $0$ is 
contained in $\cup S_i$. This implies that $f^q([p,0])\subset [p,0)$.
Then observe that $\Cal{R}_{\text{class},\sigma}(O(\Phi_k(g_k))\circ q_{t_k})$
 tends
to the rescaling of $f^q: S_q\to S_q$. So for $k$ big enough we see 
that $\Cal{R}_{\text{dyn}}(\Phi_k(g_k),t_k)$ is in $B_0$. In particular
$$
|\Cal{R}_{\text{dyn}}(\Phi_k(g_k),t_k)-(\Phi_k(g_k),t_k)|\ge \delta,
$$ 
which contradicts property 2 and 3  above.

\bigskip

\flushpar
We proved that  $(\Phi(g)), t)\in U^+$ and we can apply the dynamical 
renormalization operator. 
The continuity of the functions $\Cal{R}, d$ and $\rho$ implies
$$
\aligned
\Cal{R}_{\text{dyn}}(\Phi(g),t)&=
(\Cal{R}(\Phi(g),d(\Phi(g),t)),\rho(\Phi(g),t))\\
&=\lim_{k\to\infty}(\Cal{R}(\Phi_k(g_k), d(\Phi_k(g_k),t_k)),
                   \rho(\Phi_k(g_k),t_k))\\
&=\lim_{k\to\infty}(\Cal{R}(\Phi_k(g_k),g_k),\rho_k)\\
&=\lim_{k\to\infty}(\Phi_{k+1}(g_k),\rho_k)=(\Phi(g),t).
\endaligned
$$
In particular for $\underline{\phi}=\Phi(g)$ we get
$$
\Cal{R}_{\text{class},\sigma}(O(\underline{\phi})\circ q_t)=
O(\Cal{R}(\underline{\phi}, d(\underline{\phi},t)))
\circ q_{\rho(\underline{\phi},t)}=O(\underline{\phi})\circ q_t.
$$ 
\hfill\hfill\qed $\,\,$ (Theorem 2.2)
\enddemo

\bigskip
\centerline{\bf 8. The Existence of Truncation Fixed Points}
\bigskip

\flushpar
In this section we are going to prove Proposition 7.1, the existence of
truncation fixed points.

\bigskip

\flushpar
Assume that $(\underline{\phi},t)$ is a fixed point of $\Cal{R}_{\text{dyn}}$. 
Then
$$
\Cal{R}_{\text{dyn}}(\underline{\phi},t)=
(\Cal{R}(\underline{\phi},d(\underline{\phi},t)),\rho(\underline{\phi},t))=
(\underline{\phi},t).
$$
Hence $\underline{\phi}=\Phi(d(\underline{\phi},t))$. We have to construct a fixed 
using the collection 
$$
P=\bigcup_{\epsilon>0} \Phi(G_\epsilon)\subset \bigcup_{C<\infty} X_C.
$$ 
consisting of so-called {\it pure decompositions}. Let 
$$
U^+_P=\{(\underline{\phi},t)\in U^+| \underline{\phi}\in P\},
$$
be the collection of quasi renormalizable decomposed unimodal maps whose 
diffeomorphic part is pure. Observe that
$$
\Cal{R}_{\text{dyn}}: U^+_P\to P\times (0,\infty).
$$ 
It is this map to which we want to apply the bottom-down-top-up-Principle
discussed in the introduction. Unfortunately $P$ is homeomorphic to a countable 
product
of open intervals. It is not possible to take a compactification of $P$ and 
extend the renormalization operator to this compactification. The problem 
that arises is that the composition operator
$$
O:\bigcup_{C<\infty} X_C\to \Cal{D}
$$
can not be extended. This extension is necessary to define the dynamical 
geometry. We have to proceed differently.

\bigskip

\flushpar
Let 
$$
\Cal{Q}=\{Z_I(q_\frac12)| I=[1-a,1] \text{ an oriented interval with } a\in [0,1]\}. 
$$
This set can be parametrized by $[-1,1]$: for $a\in [0,1]$ we take the
interval $[1-a,1]$ with $o(I)=1$ and for $a\in [-1,0]$ we use $o(I)=-1$. There 
is no ambiguity when $a=0$, $Z_I(q_\frac12)=id$. 

\flushpar
This compact space $\Cal{Q}$ will be used to compactify the pure 
decompositions.
Observe that when $g\in G_\epsilon$, $\epsilon>0$,
$$
\Phi(g)(\tau)\in \Cal{Q},
$$
for all $\tau\in T$.

\bigskip

\flushpar
Instead of working in the infinite dimensional space of pure decompositions
we will work in finite dimensional truncations. 
Fix $k\ge 0$. Let
$$
P(k)=\Cal{Q}^{\cup_{j\le k} L_j}
$$
Observe that $P(k)$ is homeomorphic to a finite dimensional closed Euclidean ball.
There is an embedding $j:P(k)\to P(k+1)$ defined by
$$
\aligned
j(\underline{\phi})(\tau)&=\phi_\tau , \tau\in \cup_{j\le k} L_j\\
j(\underline{\phi})(\tau)&=id        ,\tau \in L_{k+1}
\endaligned
$$
Identify $P(k)$ with this embedding, $P(k)\subset P(k+1)$ and let
$\pi_k:P(k+1)\to P(k)$ be the projection.

\flushpar
Let $\Cal{H}$ be the collection of orientation preserving interval homeomorphisms
$h:[-1,1]\to[-1,1]$ with the property
\parindent=15pt
\item{1)} $h|_{(-1,1)}$ is $C^2$,
\item{2)} $Dh(\pm1)\ge 0$.

\flushpar
The decompositions in $P(k)$ contain only finitely many endomorphisms in $\Cal{H}$.
There is no difficulty in defining the composition operator
$$
O:P(k)\to \Cal{H}.
$$   
The definitions used before carry automatically over to the finite dimensional
truncations. The objects thus obtained can also be considered to be continuous 
extensions to the compact space $P(k)\times [0,1]$.

\flushpar
The set of quasi renormalizable decomposed unimodal maps is denoted by
$$
U^+(k)\subset P(k)\times [0,1].
$$
For $\epsilon\ge 0$ let
$$
G_\epsilon(k)=Q_\epsilon\times E_\epsilon^{\cup_{j\le k} L_j}.
$$
The dynamical geometry is again a continuous function
$$
d:U^+(k)\to \bigcup_{\epsilon>0} G_\epsilon(k)
$$
and the peak value of the renormalization is the continuous function
$$
\rho:U^+(k)\to [0,\infty).
$$
The geometrical renormalization operators 
$$
\Cal{R}_g:P(k)\to P(k+1),
$$
$g\in G_\epsilon$,  are defined as before by using the truncated set $\cup_{j\le k} L_j$
as set of decomposition times.
The dynamical renormalization operator becomes
$$
\Cal{R}_{\text{dyn}}:U^+(k)\to P(k+1)\times [0,1]
$$
with
$$
\Cal{R}_{\text{dyn}}((\underline{\phi},t))=
(\Cal{R}_{d((\underline{\phi},t))}(\underline{\phi}),\rho((\underline{\phi},t))).
$$

\bigskip

\flushpar
The proof of the following Lemma is left to the reader.

\proclaim{Lemma 8.1} For every $k\ge 0$ there exists a unique continuous 
function $\beta_k: P(k)\to (0,1]$ such that
$$
O(\underline{\phi})\circ q_{\beta_k(\underline{\phi})}(0)=0.
$$
In particular
$$
U^+(k)\subset \{(\underline{\phi},t)\in P(k)\times [0,1]| 
t> \beta_k(\underline{\phi})\}
$$
and $\beta_{k+1}|_{P(k)}=\beta_k $.
\endproclaim

\proclaim{Lemma 8.2} For $k\ge 0$
\parindent=15pt
\item{1)} $U^+(k)$ is open,
\item{2)} $P(k)\times \{1\}\subset U^+(k)$ and 
$\rho(P(k)\times \{1\})\subset (1,\infty),$
\item{3)} there is a continuous extension
$$
\Cal{R}_{\text{dyn}}:\overline{U^+(k)}\to P(k+1)\times [0,\infty).
$$ 
\item{4)} the image of the boundary of $U^+(k)$ satisfies
$$
\Cal{R}_{\text{dyn}}(\partial U^+(k))\subset 
\{(\underline{\phi},t)\in P(k+1)\times [0,1]| t<\beta_{k+1}(\underline{\phi})\}
$$
\endproclaim

\flushpar
The boundary of a set $A$ is $\partial A=\overline{A}\cap \overline{A^c}$.

\demo{Proof} Observe that a map is quasi renormalizable iff it
has the expanding periodic orbit of the combinatorial type determined by 
$\sigma$. 
The unimodal maps corresponding to points in $P(k)\times \{1\}$, 
so-called full maps, have periodic
orbits of all possible combinatorial types. All of which are expanding.
This explains Property 2).

\flushpar
The stability of expanding periodic orbits together with the above 
observation implies that $U^+(k)$ is open, property 1).  

\flushpar
The boundary of $U^+(k)$ consists of maps $(\underline{\phi},t)$ for which the
unimodal map $O(\underline{\phi})\circ q_t$ has a neutral periodic point $p$ with
combinatorics still the same as determined by $\sigma$. We can still 
construct an elementary geometry by pulling back the central interval 
$C=[-|p|,|p|]$. This elementary geometry determines the extension of 
$d$ and $\rho$. Property 3. 

\flushpar
Observe that this neutral point attracts the orbit of the critical point. The
classical renormalization of $O(\underline{\phi})\circ q_t$ is a unimodal map 
such that the image of the critical point $0$ 
is below $0$. Property 4.  
\hfill\hfill\qed $\,\,$ (Proposition 8.2)
\enddemo

\flushpar
The author would like to thank D.Sullivan for indicating the following

\proclaim{Bottum-Down-Top-Up-Proposition 8.3} 
Let $F:D_n\times [0,1]\to D_n\times \Bbb{R}$ be a
continuous map, $D_n$ is the closed $n-$dimensional ball. 
If
$$
F(D_n\times \{0\})\subset D_n\times (-\infty, 0)
$$
and
$$
F(D_n\times \{1\})\subset D_n\times (1,\infty)
$$
then $F$ has a fixed point in $D_n\times [0,1]$.
\endproclaim

\demo{Proof} Suppose not. Then let $S$ be the boundary of $D_n\times [0,1]$,
this is an $n-$dimensional sphere. And consider the displacement function
$f:S\to S^n$, where $S^n$ is the unit sphere in $\Bbb{R}^{n+1}$ and
$$
f(x)=\frac{x-F(x)}{|x-F(x)|}.
$$

\proclaim{Claim}
$deg(f)=\pm 1$.
\endproclaim

\flushpar
By radial projection onto the axis $\{0\}\times \Bbb{R}$ 
we can find a homotopy of $F$, say $F_t$, $t\in[0,1]$,
such that
\parindent=15pt
\item{1)}$F_0=F$.
\item{2)} $F_t$ has no fixed points in $S$.
\item{3)} the image of $F_1$ is an interval $\{\underline{0}\}\times [-a,b]$,
$a>0$, $b>1$.

\flushpar
Then the images $F_1(D_n\times \{0\})$ and  $F_1(D_n\times \{1\})$ are 
intervals in $\{\underline{0}\}\times [-a,b]$ outside $D_n\times [0,1]$.
Now we can perform a second homotopy to collapse these two intervals 
to points. Without loss of generality we may assume that
\item{4)} 
$$
F_1(D_n\times \{1\})=(\underline{0}, b)
$$
and
$$
F_1(D_n\times \{0\})=(\underline{0},-a).
$$

\flushpar
Let $f_t:S\to S^n$ be the displacement function for $F_t$. Because $f=f_0$ and 
$f_1$ are homotopic we get
$$
\text{deg}(f)=\text{deg}(f_1).
$$

\flushpar
To compute this degree
observe that only $(\underline{0},1)$ and  $(\underline{0},0)$ will have
pure vertical displacements under the map $F_1$. 
The first point goes up, the second goes down. 
So if $N\in S^n$ be the north pole of $S^n$ then
$$
f_1^{-1}(N)=\{(\underline{0},1)\}.
$$ 
So the degree of $f_1$ equals the local degree at $(\underline{0},1)$.
In the natural local coordinates at $(\underline{0},1)$ $f_1$ is exactly the 
anti-podal map in dimension $n$. So
$$
\text{deg}(f)=\pm 1.
$$
The Claim is proved.

\bigskip

\flushpar
But $f$ has an extension to $f:D_n\times [0,1]\to S^n$ in the obvious way,
so $\text{deg}(f)=0$. Contradiction. The map $F$ has a fixed point.
\hfill\hfill\qed $\,\,$ (Proposition 8.3)
\enddemo

\demo{The Construction of a Truncation Fixed Point, Proposition 7.1} Fix 
$k\ge 0$. Let $H_1:P(k+1)\times \Bbb{R}\to P(k+1)\times \Bbb{R}$ be defined by
$$
H_1(\underline{\phi},t)=(\underline{\phi},v(\underline{\phi},t)).
$$ 
and let $H_2:P(k+1)\times \Bbb{R}\to P(k)\times \Bbb{R}$ be defined by
$$
H_2(\underline{\phi},t)=(\pi_k(\underline{\phi}),t).
$$
Let
$$
V^+(k)=H_1(U^+(k))\subset P(k)\times (0,1].
$$ 
Consider the function $R:V^+(k)\to P(k)\times \Bbb{R}$ defined by
$$
R=H_2\circ H_1\circ \Cal{R}_{\text{dyn}}\circ H_1^{-1}.
$$

\midinsert
\centerline{\psfig{figure=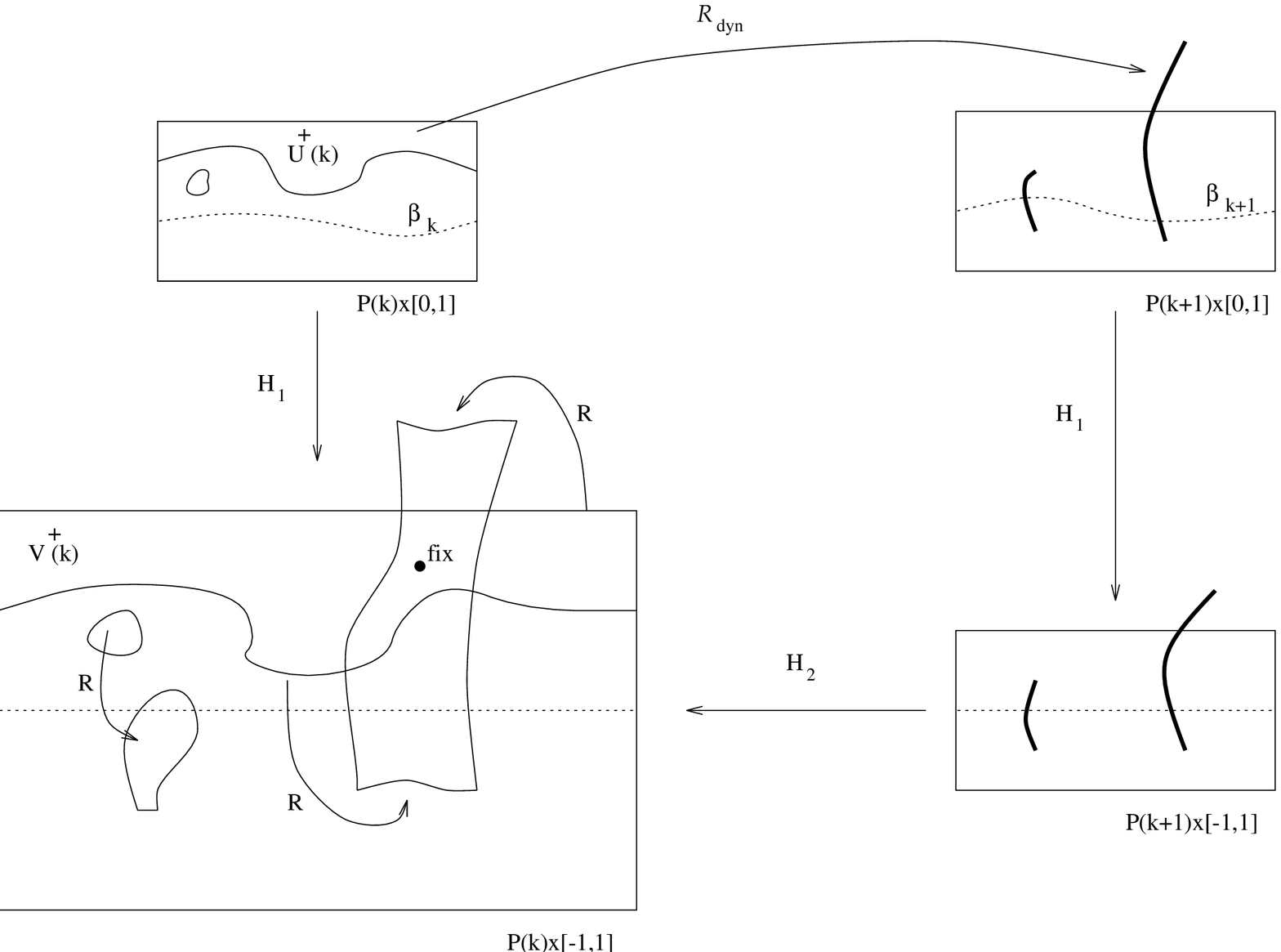,width=0.8\hsize }}
\centerline{Figure 7. $R: V^+(k)\to P(k)\times \Bbb{R}$}
\endinsert

\flushpar
It is this function to which we would like to apply the previous Proposition.
However, we do not known whether $V^+(k)$ has the right product form. It could
have different components. The solution is to extend the function $R$ 
continuously to 
$P(k)\times [0,1]$. To be able to make this extension we have to deform $R$
slightly. 

\flushpar
Let $D:P(k)\times \Bbb{R}\to P(k)\times \Bbb{R}$ defined as follows
\parindent=15pt
\item{1)} $A|_{P(k)\times [0,\infty)}=id$,
\item{2)} $A(P(k)\times (-\infty,0))\subset P(k)\times [-1,0)$,
\item{3)} $A(R(\partial V^+(k)))=\{(\underline{\phi}_0,-1)\}$.

\flushpar
The third property is possible because of Lemma 8.2(4): 
$R(\partial V^+(k))\subset P(k)\times [-1,0)$ is compact. The function
$F=A\circ R:V^+(k)\to P(k)\times \Bbb{R}$ is constant on $\partial V^+(k)$ and
hence can be continuously extended
 to $F:P(k)\times [0,1]\to P(k)\times \Bbb{R}$. 

\flushpar
It has a fixed point $(\underline{\phi},v)\in P(k)\times [0,1]$, by 
Proposition 8.3. It has to be in $V^+(k)$, otherwise it would be maps into
$P(k)\times\{-1\}$. In particular it is a fixed point for $R$. The truncation
fixed point is going to be 
$$
(\underline{\phi},t)=H^{-1}(\underline{\phi},v)\in U^+(k).
$$
There is $\epsilon>0$ such that the dynamical geometry $g$ of 
$(\underline{\phi},t)$ is in $G_\epsilon$. Because $(\underline{\phi},v)$ is a
fixed point of $R$ we have 
$$
(\underline{\phi},v)=H_2\circ H_1(\Cal{R}_{\text{dyn}}(\underline{\phi},t)).
$$
So
$$
\underline{\phi}=\pi_k(\Cal{R}_g(\underline{\phi})).
$$
Hence
$$
\underline{\phi}=\Cal{R}_g^{k+1}(id).
$$
Indeed we found a truncation fixed point in $U^+$.
\hfill\hfill\qed $\,\,$ (Proposition 7.1)
\enddemo

\bigskip
\centerline{\bf 9. A Priori Bounds}
\bigskip

\flushpar
In this section we are going to prove the a priori bounds on the 
geometry of truncation fixed points, Proposition 7.2. Fix a truncation
fixed point $(\Phi_k(g),t)\in U^+$, $k\ge 0$, say
$$
g=(\Cal{S},(\Cal{S}_\tau)_{\tau\in T})
$$
with
$$
\aligned
\Cal{S}&=\{S_1,S_2.\dots,S_q=[-|p|,|p|]\},\\
\Cal{S}_\tau&=\{S_1(\tau),S_2(\tau).\dots,S_q(\tau)\}, \text{ for } \tau\in T.
\endaligned
$$
The a priori bounds on the geometry $g$ will follow easily when a uniform
$\epsilon$ is found with
$$
\Cal{S}\in Q_\epsilon.
$$ 
The first step is to find a uniform  bounds on the position of the periodic 
point $p\in (-1,1)$ of $f=O(\Phi_k(g))\circ q_t$.

\proclaim{Lemma 9.1} 
There exist constants $K>0$ and $p_0>0$, independent of $k$
with the following property. If $|p|<p_0$ then 
$$
|g| \le K \root {\alpha\cdot q} \of {|p|},
$$
and 
$$
\sum_{\tau\in \bigcup_{j\ge \alpha q + 1} L_j} |\phi_\tau|\le K.
$$
\endproclaim

\demo{Proof} The geometry is the dynamical geometry of $(\Phi_k(g),t)$. 
It is obtained
by pulling back the interval $S_q=[-|p|,|p|]$. All the maps in $\Phi_k(g)$
have negative Schwarzian derivative which implies immediately the first 
statement. To prove the second statement observe that
$$
\sum_{i=1}^{q-1}|a_i|\le (\alpha-1)\frac1p \cdot (q-1).
$$
To recover $\Phi_k(g)$ we have to apply repeatedly the geometric 
renormalization operator $\Cal{R}_g$. The bound on
$|g|$ implies that the geometrical renormalization operator contracts with a 
constant 
$K\root {\alpha q }\of {p}$. So for all $j\ge 0$
$$
\sum_{\tau\in L_j}|\phi_\tau|\le (\alpha-1)\frac1p\cdot 
(q-1)\cdot(K\root {\alpha q} \of {p})^j.
$$
In particular
$$
\aligned
\sum_{\tau\in L_{\alpha q+1}}|\phi_\tau|&\le (\alpha-1)\frac1p\cdot 
(q-1)\cdot (K\root {\alpha q}\of {p})^{\alpha q+1}\\
&=(\alpha-1)\cdot (q-1)\cdot K^{\alpha q+1}
      \cdot p\cdot (K\root {\alpha q} \of {p})^1\le 1,
\endaligned
$$
for $p$ small enough. Because the contraction constant $K\root {\alpha q}
\of {p}$ is small, we get a uniform bound on the total norm of the 
diffeomorphisms
in levels deeper than $\alpha q+1$. 

\hfill\hfill\qed $\,\,$ (Lemma 9.1)
\enddemo

\proclaim{Lemma 9.2} There exist $K>0$ and $p_0>0$, independent of $k$, 
with the following property. If $|p|<p_0$ then
$$
(O(\Phi_k(g))'(x)\le K,
$$
for all $x\in [-1,1]$.
\endproclaim

\demo{Proof} Observe that $O(\Phi_k(g))$ is obtained by composing the maps in the 
levels $L_0,L_1,\dots, L_{\alpha q}$, which are finitely
many rescaled restrictions of the canonical folding map, 
and diffeomorphisms in the deeper levels. These diffeomorphisms in the 
deeper levels are controlled by Lemma 9.1. Unfortunately the maps in the 
first $\alpha q$ levels can be highly non-linear. However, there is a uniform 
bound on the derivative of each of them.

\flushpar
The finite number of maps in the first $\alpha q$ levels have uniform bounded 
derivative and the maps in the deeper levels are uniformly bounded, 
$O(\Phi_k(g))$ has uniform bounded derivative.
\hfill\hfill\qed $\,\,$ (Lemma 9.2)
\enddemo

\flushpar
The next Lemma states the a priori lower bound on the position of the 
periodic point $p$.

\proclaim{Lemma 9.3} There exists $\epsilon>0$, independent of $k$, such that 
$$
|p|>\epsilon.
$$
\endproclaim

\demo{Proof} Assume that the periodic point $p$ of $f$ is very close to $0$.
Because $(\Phi_k(g),t)$ is a truncation fixed point we have
$$
\left| f^q(S_q)\right| \ge \frac12 \left| S_q\right|=|p|.
$$
But
$$
\left| f^q(S_q)\right|\le (2\alpha K)^{q-1} \cdot K \cdot 2t|p|^\alpha, 
$$
where $(2\alpha K)^{q-1}$ is the bound obtained from Lemma 9.2 for the 
derivative of $f$.  These two estimates are impossible for $|p|$ very small.

\hfill\hfill\qed $\,\,$ (Lemma 9.3)
\enddemo

\proclaim{Lemma 9.4} There exists $\epsilon>0$, independent of $k$,
such that
$$
\Cal{S}\subset Q_\epsilon.
$$
\endproclaim

\demo{Proof} The previous Lemma states that the periodic point is not to 
close to $0$. Left is to show that for some uniform $\epsilon >0$
$$
\bigcup_{i=1}^q S_i \subset  (-1+\epsilon,1-\epsilon).
$$
Assume the contrary. Every interval $S_i$ lies between $S_1$ and $-S_1$.
Hence the assumption implies that the right most interval in $\Cal{S}$, 
that is $S_1$, is very close to the boundary point $1$.

\flushpar
Use the notation $\phi=O(\Phi_k(g))$ and  $\psi=O(\Phi_{k+1}(g))$ and 
$f=\phi\circ q_t$.

\bigskip

\flushpar
The decompositions $\Phi_k(g)$ and $\Phi_{k+1}(g)$ differ only in level $k+1$.
This implies that the map $\psi$ is obtained from $\phi$ by Sandwiching
the maps in the level $k+1$ of $\Phi_{k+1}(g)$ into $\phi$. From this
we can not conclude general distortion statements. However, we can compare
$\psi'(1)$ and $\phi'(1)$. They differ by the product of the derivatives in 
$1$ of the maps in this level $k+1$. 
$$
\psi'(1)=\phi'(1)\times \Pi_{\tau\in L_{k+2}} \psi'_{\tau}(1).
$$
The last factor is bounded. The reason for 
this is that the diffeomorphisms in this level $k+1$ are obtained by repeatedly
applying the geometrical renormalization operator $\Cal{R}_g$. So, their total
norm is bounded by
$$
\sum_{\tau\in L_{k+1}} |\psi_\tau| \le \sum_{i=1}^{q-1} |a_i|
\le (\alpha-1)\frac1p\cdot (q-1)\le (\alpha-1)\frac{1}{\epsilon}\cdot (q-1),
$$ 
where $\epsilon>0$ is the minimal distance of the periodic point $p$ to $0$
given in Lemma 9.3. This bound is uniform. Now apply Lemma 10.3 and we get
a uniform bound 
$$
\psi'(1)\le K\cdot \phi'(1).
$$

\flushpar
The next part compares these two derivatives dynamically. The interval
$\hat{S}_1=\phi^{-1}(S_1)$ has two boundary points: $q_t(p)$ and $b_0$. Let
$b_1=\phi(b_0)$ be the corresponding boundary point of $S_1$. Then
$$
\psi'(1)=\frac{|\hat{S}_1|}{|S_q|}\cdot (f^{q-1})'(b_1)\cdot \phi'(b_0). 
$$
The first factor is a normalization factor. Compare the dynamical picture
in Figure 6. We are going to show a uniform bound
$$
(f^{q-1})'(b_1)\le B.
$$
This will be done in a few steps.

\bigskip

\flushpar
First we claim 
$$
|S_1|\approx 0.
$$ 
If $S_1$ is big then the periodic point $p$ has to be away from the boundary 
and the interval $S_q=[-|p|,|p|]$ has bounded hyperbolic length. The other 
intervals $S_i$ are obtained by pulling back by $\phi$ and $q_t$. This maps
has negative Schwarzian derivative and we get a uniform bound on the 
hyperbolic length
of all intervals $S_i$, $i=1,2,\dots, q$. But $S_1$ is very close to the 
boundary (by assumption) and $|S_1|$ is big. Contradiction.

\bigskip

\flushpar
Next we claim a uniform lower bound
$$
|\hat{S}_1|\ge \delta.
$$
The critical value of $(\Phi_k(g),t)$ is denoted by $v\in [-1,1]$. Also
use the notation $\rho=\rho(\Phi_k(g),t)$. Because $(\Phi_k(g),t)$ is a 
truncation fixed point we have 
$$
2t-1=\phi^{-1}(v) \text{ and } 2\rho-1=\psi^{-1}(v).
$$
From Lemma 10.6 (see appendix) we get a universal constant $K$ 
such that the hyperbolic distance between $2t-1\in [-1,1]$ and 
$2\rho-1\in[-1,1]$ is bounded by
$$
K\sum_{\tau\in L_{k+1} }|\phi_\tau|\le 
K\sum_{i=1}^{q-1} |a_i|\le K\cdot 
(\alpha-1)\cdot \frac{1}{\epsilon} \cdot (q-1),
$$
where $\epsilon>0$ is given by Lemma 9.3. We get universal bound on the 
hyperbolic distance between $t,\rho\in [0,1]$. In particular if $t$ is very 
small we get 
$$
\rho\le \text{const}\cdot  t. 
$$
Now we can finish the lower estimate for 
$$
\left| \hat{S}_1\right|=\frac{|q_t(S_q)|}{\rho}
           =\frac{2t\cdot p^{\alpha}}{\rho}\ge \frac{t}{\rho} \cdot 
2\epsilon^\alpha,
$$
where $\epsilon>0$ is given again by Lemma 9.3.
In case $t$ is not too small there will be a lower bound because $\rho\le 1$.
For small $t$ we use the universal estimate $\rho\le \text{const}\cdot t$ to 
obtain a universal lower bound for $|\hat{S_1}|$.

\bigskip

\flushpar
At last we claim 
$$
\phi'(b_0)\ge \phi'(1).
$$ 
The interval $q_t(S_1)$ is very small. Because $S_q$ is not small there has to
be a point $x$ between $-1$ and $\hat{S}_1$ with $\phi'(x)>>1$. 
The Minimal Principle for maps with negative Schwarzian derivative 
(see [MS]) applied to $\phi :[x,b_0]\to[-1,1]$ implies that the average
slope on $\hat{S}_1$ satisfies
$$
\phi'(x)>> 1 >> \frac{S_1}{\hat{S}_1}\ge \min\{\phi'(x),\phi'(b_0)\}.
$$
So
$$
\phi'(x)>\phi'(b_0).
$$
The minimal Principle applied to $\phi :[x,1]\to[-1,1]$ implies
$$
\phi'(x)>\phi'(b_0)\ge \min\{\phi'(x),\phi'(1) \}.
$$ 
The claim follows.

\bigskip

\flushpar
These estimates together give
$$
K\phi'(1)\ge \psi'(1)\ge {\text const}\cdot (f^{q-1})'(b_1)\cdot \phi'(1).
$$
Indeed we get a uniform bound for $(f^{q-1})'(b_1)$.

\bigskip

\flushpar
The third part of the proof will produce the contradiction. Let 
$H_i\subset [-1,1]$, $i=1,2\dots, q-1$, be the connected component of 
$[-1,1]\setminus S_i$ not containing $0$. We claim
$$
\frac{|H_i|}{|S_i|}<<1
$$
for all $i=1,2,\dots,q-1$. Assume the contrary: there is an $S_i$ with bounded
hyperbolic length. This interval can be pulled back to show that also 
$\hat{S}_1$ has bounded hyperbolic length. One step more shows that also
$S_q$ has bounded hyperbolic length. Continue to pull back and it turns out 
that every interval $S_j$ has bounded hyperbolic length. This implies that
$$
f^{q-1}:S_1\to S_q
$$
has bounded distortion. Because $\frac{|S_q|}{|S_1|}$ is very big we get that
$(f^{q-1})'(b_1)$ is also very big. Contradiction.

\bigskip

\flushpar
This Claim implies that for every $i=2,3,\dots,q-1$ we have 
$S_1\subset H_i\cup -H_i$ and 
$$
|S_1|\le |H_i|<< |S_i|.
$$

\bigskip

\flushpar
Let $A=(b_1,a]$ be the maximal interval on which $f^{q-1}$ is monotone.
The map $f^{q-1}$ has a very big average slope on $S_1$ but
the derivative in $b_1$ is uniformly bounded by $B$. The Minimal Principle
implies that 
$$
(f^{q-1})'(z)\le B,
$$  
for every point $z\in A$. We will construct a point in $A$ with a very big 
derivative and so produce a contradiction.

\bigskip

\flushpar
Observe
$$
|f^{q-1}(A)|\le B\cdot |A|<< |S_1|,
$$
$f^{q-1}(A)$ does not contain $S_1$ (or $-S_1$). This cannot happen in the
period doubling case $q=2$. So $q\ge 3$. Moreover, there has to be some 
$j\le q-2$ with
$$
f^j(a)=0.
$$

\midinsert
\centerline{\psfig{figure=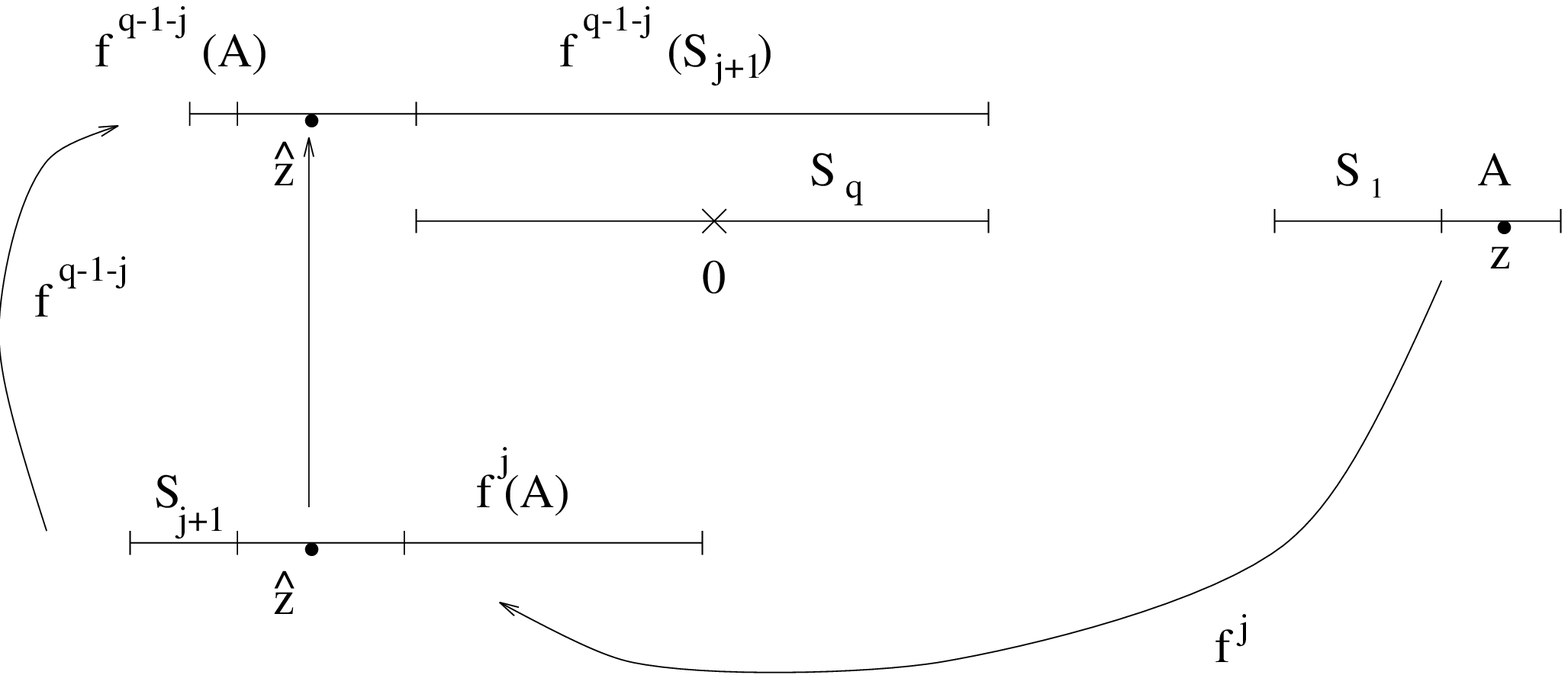,width=0.6\hsize }}
\centerline{Figure 8.}
\endinsert

\flushpar
There has to be a point $z\in A$ such that $\hat{z}=f^j(z)$ (or $-\hat{z}$)
is periodic with period $q-1-j< q$, see Figure 8. Such a point can not be in 
the orbit of 
$S_1$ because this orbit has period $q$. In particular it does not attract the critical orbit. The map $f$ has negative Schwarzian derivative which implies that 
the periodic orbit is expanding
$$
|(f^{q-1-j})'(\hat{z})|>1.
$$
Let $A'=f^{-j}(S_q)\cap A$ and estimate
$$
\aligned
|(f^{q-1})'(z)|&=|(f^{q-1-j})'(\hat{z})|\cdot |(f^{j})'(z)|\\
               &\ge 1\cdot |(f^{j})'(z)|\\
               &\ge \min\{\frac{|f^j(S_1)|}{|S_1|},
                          \frac12\cdot\frac{|S_q|}{|A'|}\}\\
               &\ge \min\{\frac{|S_{j+1}|}{|S_1|},
                          \frac12\cdot\frac{|S_q|}{|S_1|}\}>>1.
\endaligned
$$
Contradiction.
\hfill\hfill\qed $\,\,$ (Lemma 9.4)
\enddemo

\demo{Proof of the A Priori Bounds, Proposition 7.2} The previous Lemma gives 
uniform bound on the elementary geometry $\Cal{S}\in Q_\epsilon$. It is 
contained in an interval with bounded hyperbolic length. The 
Schwarzian derivative of all the diffeomorphisms $O^\tau(\Phi_k(g))$ is 
negative we can pull back this elementary geometry $\Cal{S}$ and recover the
whole geometry $g$ and see that all the elementary geometries $\Cal{S}_\tau$
are contained in an interval with uniformly bounded hyperbolic length. 
In particular they have also Euclidean length uniformly bounded away from 2,
$$
|g|\le 1-\delta.
$$ 
The geometric renormalization operator $\Cal{R}_g$ has 
uniformly bounded contraction constant. This implies that the decomposition
$\Phi_k(g)$ lies in a uniformly bounded set in $X_C$. We get a uniform bound
on the derivatives of the diffeomorphisms $O^\tau(\Phi_k(g))$. Hence 
the dynamical
intervals $S_i(\tau)$ with $i=1,2,\dots q$ and $\tau \in T$ can also not be 
to close to the boundary of $[-1,1]$. The a priori bounds on the geometry of 
any truncation fixed point are shown.

\hfill\hfill\qed $\,\,$ (Proposition 7.2)
\enddemo

\bigskip
\centerline{\bf 10. Appendix}
\bigskip

\flushpar
In this appendix we discuss some Lemmas describing the relation between
the $C^2$ topology on $\Cal{D}=\text{Diff}_+^2([-1,1])$ and the non-linearity norm used
throughout the text.

\proclaim{Lemma 10.1(Chain-Rule for Non-linearities)} Let $\phi,\psi\in 
\Cal{D}$ then
$$
\eta_{\psi\circ \phi}(x)=\eta_{\psi}(\phi(x))\cdot \phi'(x)+
\eta_\phi(x),
$$
for all $x\in [-1,1]$.
\endproclaim

\proclaim{Lemma 10.2} If $\phi,\psi\in \Cal{D}$ then
$$
\left| \psi(x)-\phi(x)\right| \le 2(e^{4|\psi-\phi|}-1)
$$
for all $x\in [-1,1]$.
\endproclaim

\demo{Proof} Use the inverse of the non-linearity:
$$
\aligned
\psi(x)&=2\frac{\int_{-1}^x e^{\int_{-1}^s\eta_\psi}ds}
                      {\int_{-1}^1 e^{\int_{-1}^s\eta_\psi}ds}-1\\
       &=2\frac{\int_{-1}^x e^{\int_{-1}^s(\eta_\psi-\eta_\phi+\eta_\phi)}ds}
                      {\int_{-1}^1 e^{\int_{-1}^s
         (\eta_\psi-\eta_\phi+\eta_\phi)}ds}-1\\
       &\le e^{4|\psi-\phi|}\cdot \phi(x) +e^{4|\psi-\phi|}-1.
\endaligned
$$
So
$$
\aligned
\psi(x)-\phi(x) &\le (e^{4|\psi-\phi|}-1)\cdot \phi(x) +e^{4|\psi-\phi|}-1\\
                &\le 2(e^{4|\psi-\phi|}-1).
\endaligned
$$
The above estimate is symmetric in $\psi$ and $\phi$: the Lemma follows.
\hfill\hfill\qed $\,\,$ (Lemma 10.2)
\enddemo

  \proclaim{Lemma 10.3} Let $\psi$ be a composition of finitely many 
  $\phi_i\in\Cal{D}$, $i=1,2,\dots, s$. Let 
  $|\underline{\phi}|=\sum |\phi_i|$. Then
  $$
  e^{-2\left|\underline{\phi}\right|}\le 
  \psi'(x)\le e^{2\left|\underline{\phi}\right|}
  $$
  and
  $$
  e^{-
  \left|\underline{\phi}\right|
	      e^{2\left|\underline{\phi}\right|}\cdot \left|x-y\right|}
  \le
  \frac{\psi'(x)}{\psi'(y)}
  \le 
  e^{
  \left|\underline{\phi}\right|
	      e^{2\left|\underline{\phi}\right|}\cdot \left|x-y\right|}
  $$
  for all $x,y\in[-1,1]$. 
  \endproclaim

  \demo{Proof} Let $\psi=\phi_s\circ\dots\circ\phi_2\circ\phi_1$. Take
  $x,y\in [-1,1]$ and let $x_i$ and $y_i$ be the images of $x$ and $y$ under
  the partial composition $\phi_{i-1}\circ\dots\circ\phi_2\circ\phi_1$. Then
  $$
  \aligned
  \left|\ln\psi'(x)-\ln\psi'(y)\right|
  &=\left|\sum_{i}\ln\phi'_i(x_i)-\ln\phi'_i(y_i)\right|\\
  &\le \sum_{i}\left| \ln\phi'_i(x_i)-\ln\phi'_i(y_i)\right|\\
  &\le \sum_i \left|\phi_i\right| |x_i-y_i|\\
  &\le \left|\underline{\phi}\right|\cdot 2.
  \endaligned
  $$ 
  Because there is some point where $\psi'$ equals $1$ we get the first 
  estimate. Moreover we can use this estimate to get a bound 
  $$
  |x_i-y_i|\le e^{2\left|\underline{\phi}\right|}|x-y|.
  $$
  This gives us immediately the second estimate.
  \hfill\hfill\qed $\,\,$ (Lemma 10.3)
  \enddemo

\proclaim{Lemma 10.4} For every bounded set $B\subset \Cal{D}$ there exists a 
constant $K$ such that for any pair $\psi,\phi\in B$
$$
|y-x|\le K\cdot |\psi(y)-\phi(x)|+ K\cdot |\psi-\phi|
$$
for all $x,y\in [-1,1]$.
\endproclaim

\demo{Proof}
$$
\aligned
\psi(y)&=\psi(x)+\psi'(\theta)(y-x)\\
       &=\phi(x)+\psi'(\theta)(y-x)+\psi(x)-\phi(x)\\
\endaligned
$$
So
$$
\aligned
|y-x| &=\left| \frac{1}{\psi'(\theta)}\cdot 
                      \{ \psi(y)-\phi(x)+(\phi(x)-\psi(x))\}\right|\\
    &\le K\cdot |\psi(y)-\phi(x)|+K\cdot |\phi-\psi| 
\endaligned
$$
where we used Lemma 10.2 and Lemma 10.3.
\hfill\hfill\qed $\,\,$ (Lemma 10.4)
\enddemo

  \proclaim{The Sandwich Lemma 10.5} For all $b>0$ and $C>0$ 
there exists a Sandwich 
  constant $K$, 
  such that the following holds. Let $\psi_1,\psi_2$ be compositions of
  finitely many $\phi,\phi_i\in \Cal{D}$, $i=1,2,\dots,s$:
  $$
  \psi_1=\phi_s\circ\dots\circ\phi_t\circ\phi_{t-1}\circ\dots
	 \phi_2\circ\phi_1
  $$
  and 
  $$
  \psi_2=\phi_s\circ\dots\circ\phi_t\circ\phi\circ\phi_{t-1}\circ\dots
	 \phi_2\circ\phi_1.
  $$
  If $\sum_{i} |\phi_i|+|\phi|\le b$ and for $i=1,\dots,s$ 
  $|\eta'_{\phi_i}(x)|\le C |\eta_{\phi_i}(x)|$ then
  $$
  |\psi_2-\psi_1|\le K |\phi|.
  $$
  \endproclaim

  \demo{Proof} Let $x\in [-1,1]$. For $0\le i\le s-1$ define 
  $$
  x_i=\phi_{i}\circ\dots\circ\phi_2\circ\phi_1(x) \text{ ,  } x_0=x,
  $$ 
  and 
  $$
  D_i=(\phi_{i}\circ\dots\circ\phi_2\circ\phi_1)'(x).
  $$  
  Furthermore for $t-1\le j\le s-1$ let
  $$
  x'_j=\phi_j\circ\dots\circ\phi_t\circ\phi(x_{t-1})
  $$ 
  and 
  $$
  D'_j=(\phi_{j}\circ\dots\circ\phi_t)'(x'_{t-1})\phi'(x_{t-1})D_{t-1}.
  $$
  For $0\le j\le t-2$ let
  $$
  x'_j=x_j \text{ and } D'_j=D_j.
  $$

  \flushpar
  Then the chain rule for non-linearities gives
  $$
  \aligned
  |\eta_{\psi_2}(x)-\eta_{\psi_1}(x)|
  &=\left|\sum_{i=0}^{s-1} \{
  \eta_{\phi_{i+1}}(x'_i)D'_i-\eta_{\phi_{i+1}}(x_i)D_i\}+
  \eta_{\phi}(x_{t-1})D_{t-1} \right|\\
  &\le \sum_{i=t-1}^{s-1} \left|
  \eta_{\phi_{i+1}}(x'_i)D'_i-\eta_{\phi_{i+1}}(x_i)D_i\right|+
  \left|\eta_{\phi}(x_{t-1})D_{t-1}\right|.
  \endaligned
  $$
  The last term is easy to estimate. From Lemma 10.3 we have
  a uniform estimate on the derivatives $D_i$,
$$
  \left|\eta_{\phi}(x_{t-1})D_{t-1}\right|\le K\cdot |\phi|.
$$ 
Let us concentrate on the other 
  terms. We will use the symbol $K$ for all constants appearing in the 
estimates.
  It will only depend on $b$ and $C$.

  \flushpar
  We need two other derivatives: let $t-1\le i\le s-1$ then
  $$
  E_i=(\phi_{i}\circ\dots\circ\phi_t)'(x_{t-1})
  $$
  and 
  $$
  E'_i=(\phi_{i}\circ\dots\circ\phi_t)'(x'_{t-1}).
  $$
  Then 
  $$
  D_i=E_i\cdot D_{t-1} \text{ and } D'_i=E'_i\cdot \phi'(x_{t-1})\cdot D_{t-1}
  $$
  for $t-1\le 1\le s-1$.
  Take $t-1\le i\le s-1$ and consider the corresponding term
  $$
  \left|\eta_{\phi_{i+1}}(x'_i)D'_i -\eta_{\phi_{i+1}}(x_i)D_i\right|=
  $$
  $$
  \left|\eta_{\phi_{i+1}}(x'_i)E'_i\phi'(x_{t-1}) -
	  \eta_{\phi_{i+1}}(x_i)E_i\right|D_{t-1}=
  $$
  $$
  \left| (\eta_{\phi_{i+1}}(x_i)+\eta'_{\phi_{i+1}}(\theta)(x'_i-x_i))
	   \frac{E'_i}{E_i}E_i\phi'(x_{t-1})-
	    \eta_{\phi_{i+1}}(x_i)E_i\right|D_{t-1}=
  $$
  $$
  \left| \eta_{\phi_{i+1}}(x_i)
	 \{\frac{E'_i}{E_i}\phi'(x_{t-1})-1    \}E_i +
	     \eta'_{\phi_{i+1}}(\theta)(x'_i-x_i)
		  E'_i\phi'(x_{t-1}) \right|D_{t-1}\le
  $$
  $$
  |\phi_{i+1}| \cdot \left| \frac{E'_i}{E_i}\phi'(x_{t-1})-1\right| \cdot K+
  C\cdot |\phi_{i+1}|\cdot K \cdot \phi'(x_{t-1}) \cdot |x'_i-x_i|.
  $$

  \flushpar
  To continue we have to estimate $|x'_i-x_i|$. The composition 
  $\phi_{i-1}\circ\dots\circ\phi_t$ has by Lemma 10.3 a derivative bounded by
  $K$. So $|x'_i-x_i|\le K \cdot |x'_{t-1}-x_{t-1}|$. By using Lemma 10.2 we
can estimated this distance in terms of the norm of $\phi$: 
 This gives us a constant such that 
$|x'_{t-1}-x_{t-1}|\le K \cdot |\phi|$.
Hence we get a constant such that
$$
C\cdot |\phi_{i+1}|\cdot K \cdot \phi'(x_{t-1}) \cdot |x'_i-x_i| 
\le K \cdot |\phi_{i+1}|\cdot |\phi| .
$$
It is left to estimate $\left| \frac{E'_i}{E_i}\phi'(x_{t-1})-1\right|$. By 
Lemma 10.3 
$$
e^{-K|x'_{t-1}-x_{t-1}|}\cdot e^{-2|\phi|}-1\le
\frac{E'_i}{E_i}\phi'(x_{t-1})-1\le
e^{K|x'_{t-1}-x_{t-1}|}\cdot e^{2|\phi|}-1.
$$
Hence
$$
e^{-K |\phi|}\cdot e^{-2|\phi|}-1\le
\frac{E'_i}{E_i}\phi'(x_{t-1})-1\le
e^{K |\phi|}\cdot e^{2|\phi|}-1.
$$
The factor $\left|\frac{E'_i}{E_i}\phi'(x_{t-1})-1\right|$ can be 
estimated by a constant times the norm of $\phi$.
Taking all estimates together we get
$$
\left|\eta_{\phi_{i+1}}(x'_i)D'_i -\eta_{\phi_{i+1}}(x_i)D_i\right|\le
K\cdot |\phi| \cdot |\phi_{i+1}|.
$$
Moreover
$$
|\psi_2-\psi_1|=\sup |\eta_{\psi_2}(x)-\eta_{\psi_1}(x)|\le K \cdot |\phi|.
$$
The Sandwich Lemma is proved.
\hfill\hfill\qed $\,\,$ (Lemma 10.5)
\enddemo

\proclaim{Lemma 10.6} Let $\psi_1,\psi_2$ be compositions of
  finitely many $\phi,\phi_i\in \Cal{D}$ which expand hyperbolic distance, 
$i=1,2,\dots,s$:
  $$
  \psi_1=\phi_s\circ\dots\circ\phi_t\circ\phi_{t-1}\circ\dots
	 \phi_2\circ\phi_1
  $$
  and 
  $$
  \psi_2=\phi_s\circ\dots\circ\phi_t\circ\phi\circ\phi_{t-1}\circ\dots
	 \phi_2\circ\phi_1.
  $$
There is a universal constant $K$ such that the hyperbolic distance 
between
$\psi_1^{-1}(v)$ and $\psi_2^{-1}(v)$, $v\in (-1,1)$ is bounded by 
$K|\phi|$.
\endproclaim

\demo{Proof} The Lemma follows immediately once it is proved that the
hyperbolic distance between  $\phi^{-1}(v)$ and $v$ is proven to be bounded
by $K|\phi|$. From Lemma 10.3 we get $e^{-2|\phi|}\le \phi'(x)\le e^{2|\phi|}$.
Then we can use the figure 9. to estimate the maximal distance between
$v$ and $\phi^{-1}(v)$. A computation finishes the proof of the Lemma.
\hfill\hfill\qed $\,\,$ (Lemma 10.6)
\enddemo

\midinsert
\centerline{\psfig{figure=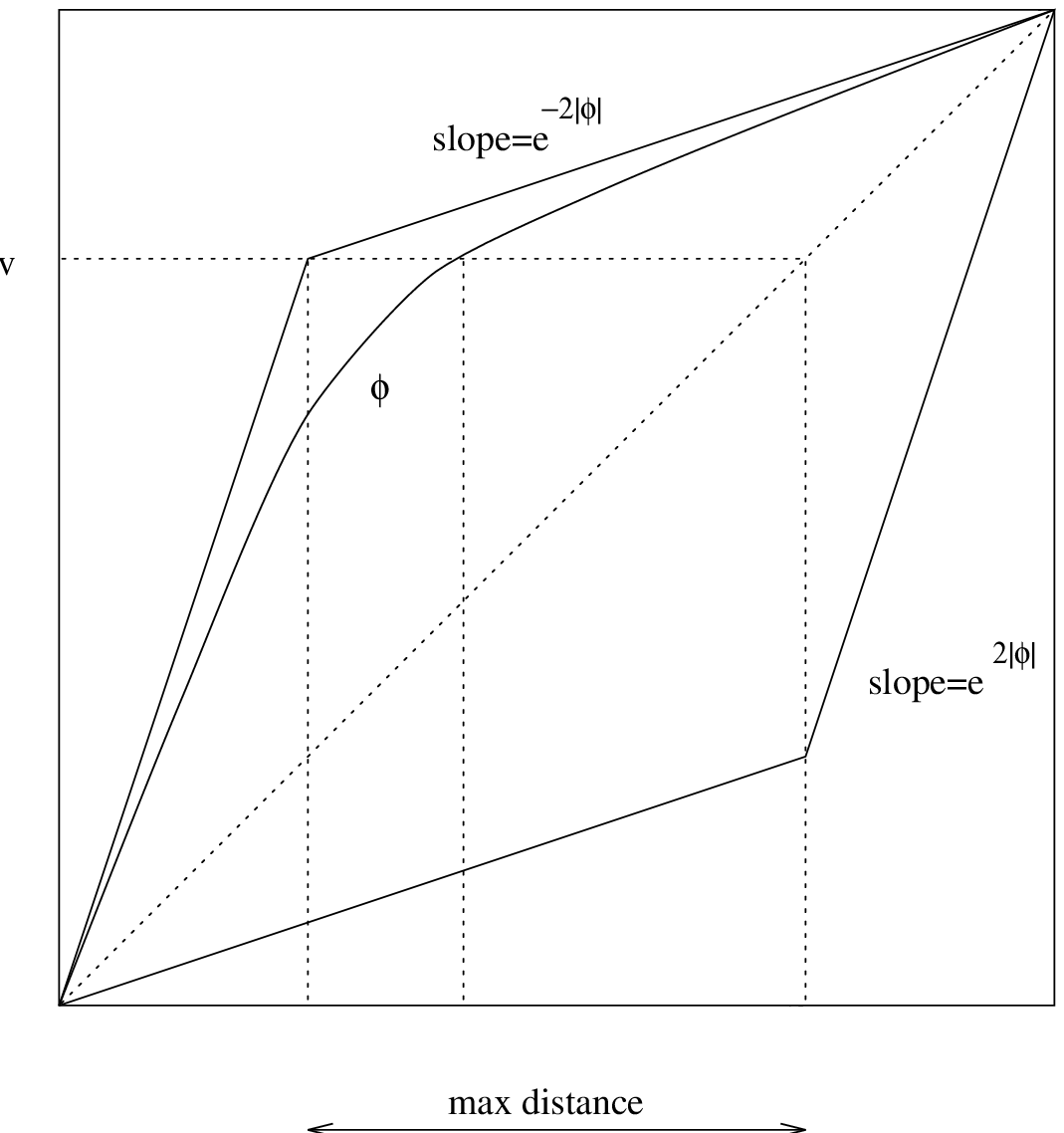,width=0.6\hsize }}
\centerline{Figure 9.}
\endinsert

\bigskip
\centerline{\bf References}
\bigskip

\parindent=30pt
\item{[CT]} P.Coullet and C.Tresser, {\it Iteration d'endomorphismes et 
            groupe de renormalisation}, J.Phys.Colloq. C {\bf 539}, C5-25 
            (1978), C.R. Acad. Sci. Paris {\bf 287} A (1978).
\item{[E1]} H.Epstein, 
           {\it New Proofs of the existence of the Feigenbaum functions},
           Comm. Math. Phys., {\bf 106} (1986), 395-426. 
\item{[E1]} H.Epstein,
           {\it Fixed points of composition operators II},
           Non-linearity {\bf 2} (1989) 305-310.
\item{[F]}  M.J.Feigenbaum, {\it Qualitative universality for a class of 
           non-linear transformations}, J. Stat. Phys. {\bf 21}, 669-706.
\item{[L]} O.E. Lanford III, {\it A computer-assisted proof of the Feigenbaum
           conjectures}, Bull.Amer.Math.Soc. New Series {\bf 6}, 1984.
\item{[Ly]} M.Lyubich, {\it Feigenbaum-Coullet-Tresser Universality and 
           Milnor's Hairiness Conjecture}, preprint.
\item{[M]} C.McMullen, {\it Complex Dynamics and Renormalization},  
           Ann. of Math., Studies {\bf 135}, Princeton University Press.
\item{[MS]} W.de Melo and S.van Strien, {\it One-dimensional Dynamics},
           Springer-Verlag, 1993.
\item{[S]} D.Sullivan, {\it Bounds, quadratic differentials,  and 
           renormalization conjectures}, Centennial Symposium, 
           Amer.Math.Soc. 1992, 417-466.

\bye